\documentclass[a4paper,leqno,11pt]{amsart}

\usepackage{ebgaramond}

\usepackage{setspace}
\usepackage{caption}
\usepackage{csquotes}
\usepackage{amsmath,amsthm}
\usepackage{amssymb}
\usepackage{xspace}
\usepackage{bm}
\usepackage{amstext}
\usepackage{comment}
\usepackage{amsfonts}
\usepackage{graphicx}
\usepackage[mathscr]{euscript}
\usepackage{amscd}
\usepackage{latexsym}
\usepackage{amssymb}
\usepackage{enumerate}
\usepackage{mathrsfs}
\usepackage[cmtip,all]{xy}
\usepackage{xcolor}
\usepackage{standalone}
\usepackage{tikz}
\usetikzlibrary{arrows.meta}
\usepackage{tikz-cd}
\usetikzlibrary{intersections}
\usetikzlibrary{calc}
\usetikzlibrary{decorations.markings}
\usetikzlibrary{arrows}
\usetikzlibrary{shapes,arrows,shapes.multipart}

\usepackage{xcolor}
\usepackage{enumitem}

\usepackage[colorlinks=true,linkcolor=blue,citecolor=blue]{hyperref}

\usepackage[cmtip,all]{xy}

\newtheorem{thm}{Theorem}[section]
\newtheorem{prop}[thm]{Proposition}
\newtheorem{lemma}[thm]{Lemma}
\newtheorem{cor}[thm]{Corollary}

\theoremstyle{definition}
\newtheorem{defn}[thm]{Definition}
\newtheorem{example}[thm]{Example}
\newtheorem{remark}[thm]{Remark}

\newtheorem*{thma}{Theorem A}
\newtheorem*{cora}{Corollary A.1}
\newtheorem*{thmb}{Theorem B}
\newtheorem*{thmc}{Theorem C}
\newtheorem*{thmd}{Theorem D}
\newtheorem*{thme}{Theorem E}
\newtheorem*{thmf}{Theorem F}

\newcommand{\uZ}{\underline{\mathbb{Z}}}
\newcommand{\uM}{\underline{M}}
\newcommand{\uN}{\underline{N}}
\newcommand{\bF}{\mathbb{F}_2}
\newcommand{\uF}{\underline{\mathbb{F}}_2}

\newcommand{\ug}{\underline{g}}
\newcommand{\gDR}{\underline{g}_{DR}}

\newcommand{\upi}{\underline{\pi}}
\newcommand{\Sp}{\mathsf{Sp}}
\newcommand{\Burn}{\mathsf{Burn}}
\newcommand{\Mor}{\mathrm{Mor}}
\newcommand{\op}{\mathrm{op}}

\newcommand{\res}{\mathrm{res}}
\newcommand{\tr}{\mathrm{tr}}
\newcommand{\aO}{\alpha_0}
\newcommand{\aI}{\alpha_1}
\newcommand{\aII}{\alpha_2}
\newcommand{\EFt}{\widetilde{E\mathcal{F}}}

\newcommand{\F}{\bF}
\newcommand{\Ab}{\mathsf{Ab}}
\newcommand{\cF}{\mathcal{F}}

\newenvironment{mysubsection}[1]
{%
  \refstepcounter{thm}%
  \par\medskip
  \noindent\textbf{\thethm.\ #1.}\upshape\ignorespaces
}
{\par\medskip}

\newenvironment{myeq}
{%
  \refstepcounter{thm}%
  \begin{equation}
  \tag{\thethm}
}
{%
  \end{equation}
}

\newcommand{\BZ}{\mathcal{B}^{\mathbb Z}}
\newcommand{\ZMod}{\uZ\text{-}\mathrm{Mod}}
\newcommand{\Lan}{\mathrm{Lan}}

\newcommand{\Hb}{\overline{H}}
\newcommand{\pr}[3]{\operatorname{pr}_{[#2,\,#3]}\big(#1\big)}   

\newcommand{\uf}{\underline{f}}

\newcommand{\im}{\operatorname{im}}
\newcommand{\Hom}{\operatorname{Hom}}
\newcommand{\coker}{\operatorname{coker}}

\begin{document}
\title{Klein four slices of the representation suspensions of $H\underline{\mathbb{F}_2}$}

\author{Surojit Ghosh}
\address{Department of Mathematics, Indian Institute of Technology, Roorkee, Uttarakhand-247667, India}
\email{surojit.ghosh@ma.iitr.ac.in; surojitghosh89@gmail.com}

\author{Ankit Kumar}
\address{Department of Mathematics, Indian Institute of Technology, Roorkee, Uttarakhand-247667, India}
\email{ankit\_k@ma.iitr.ac.in}

\date{\today}

\begin{abstract}
We study the regular slice filtration of representation spheres $S^V\wedge H\underline{M}$. For every finite group $G$, we
locate all slices in an explicit interval $[n(V),k(V)]$ and show that both endpoints occur. We introduce a \emph{defect vector} to measure a slice's failure.  For the Klein four-group with $\underline{\mathbb F}_2$-coefficients, a balanced defect allows us to move
$S^V\wedge H\underline{\mathbb F}_2$ by certain $RO(G)$-graded cohomology classes to its $\dim V$-slice, which is again a representation sphere, while an odd defect obstructs any spherical slice. As an application, we obtain closed formulas for the bottom slice of $\Sigma^nH\underline{\mathbb{F}}_2$ for every $n$, recovering the computation of Guillou and Yarnall by
a different method. Along the way, we determine the complete slice tower of $S^{\infty\alpha_i}\wedge H\underline{\mathbb F}_2$ and consequently, the partial tower for $S^{n\alpha_i}\wedge H\underline{\mathbb F}_2$.
\end{abstract}

\maketitle

\section{Introduction}
The slice spectral sequence has proved to be one of the most effective tools for computing the equivariant homotopy groups of genuine equivariant spectra. Its first major success was solving the Kervaire invariant one problem \cite{HHR16}, which uses the computation of slices of the norm $N_{C_2}^{C_8}MU_{\mathbb R}$, together with the vanishing of certain slice groups in a range. Since then it has been used to compute fixed points of Lubin--Tate theories with group actions \cite{HSWX23, BHSZ21}, to determine Hurewicz images of Real bordism \cite{LSWX19}, to relate norms of Real bordism to the Segal conjecture \cite{MSZ20}, and to analyse quotients of these norms \cite{BHLSZ22}. For cyclic groups, slice towers
have been computed in \cite{HHR16,HHR17b,HY18,Yar17,Zen17,Gho22}. Beyond cyclic groups, slices have been computed for quaternionic Eilenberg--Mac Lane spectra and quaternionic Real bordism \cite{GS25,GKM26}, and the $RO(K_4)$-graded homotopy of Klein-four normed Mackey functors has been determined \cite{GKM25}. In all of these computations, the $E_2$-page is built from the $RO(G)$-graded homotopy of the slices. So progress depends on knowing the slice towers of the simplest equivariant spectra, the representation spheres $X(V,\uM):=S^V\wedge H\uM$, as $V$ ranges over all of $RO(G)$.  For the Klein four-group, $K_4$, the smallest non-cyclic group and the first case beyond the cyclic theory, this input has been missing except for integer suspensions \cite{GY20, Slo22}.

\subsection{The slice interval}~
Our first result holds for every finite group. It locates all slices of $S^V\wedge H\uM$ by a formula in the fixed-point dimensions of $V$ and two families of integers attached to $\uM$. For $H\le G$ let
\[
b_H(\uM)=\min\{q:\pi_q\Phi^HH\uM\neq0\}\quad\text{  and  }\quad c_H(\uM)=b_H(\uM^\vee)
\]
where $\Phi^H$ denotes geometric fixed points, $\uM^\vee=\operatorname{Hom}(\uM,\mathbb Q/\mathbb Z)$ is the dual Mackey functor (restriction and transfer interchanged), and $b_H=\infty$ if $\Phi^HH\uM\simeq*$. Write $d_H(V)=\dim V^H$.
 
\begin{thma}\label{tA}
(a) Let $G$ be a finite group, $\uM$ a nonzero $G$-Mackey functor and $V\in RO(G)$. Then all slices of $S^V\wedge H\uM$ lie in levels $[n(V),k(V)]$, where
\[
n(V)=\min_{H\le G}|H|\big(d_H(V)+b_H(\uM)\big),\qquad k(V)=\max_{H\le G}|H|\big(d_H(V)-c_H(\uM)\big),
\]
 with $P^{n(V)}_{n(V)}(S^V\wedge H\uM)\neq*$ and $P^{k(V)}_{k(V)}(S^V\wedge H\uM)\neq*$. If $\uM(G/e)\neq0$, then $n(V)\le\dim V\le k(V)$, and $S^V\wedge H\uM$ is a $\dim V$-slice if and only if $n(V)=k(V)$.

(b) Moreover, if $\uM$ be cohomological, then 
\[
n(V)=\min_{H \in \mathcal{P}(G)}  |H|\big(d_H(V)+b_H(\uM)\big),\qquad k(V)=\max_{H \in \mathcal{P}(G)}|H|\big(d_H(V)-c_H(\uM)\big),
\]
where $\mathcal{P}(G)$ denote the set of subgroups of $G$ of prime-power order, including the trivial subgroup. 
\end{thma}\label{cA}
See Theorem~\ref{trange} and Corollary~\ref{ccohom}.

For the Klein four-group, $K_4:=\langle a,b|a^2=b^2=1\rangle$, we denote its three order-two subgroups by
$H_0=\langle a\rangle$, $H_1=\langle ab\rangle$, and $H_2=\langle b\rangle$.
Let $\alpha_i$ denote the nontrivial irreducible representation with kernel
$H_i$. For the constant Mackey functor $\uM=\uF$ and $G=K_4,$ we obtain the following corollary.
\begin{cora}
For $V\in RO(K_4)$ with $n=\dim V$,
\begin{align*}
&n(V)=\min\big(n,\,2d_{H_0}(V),\,2d_{H_1}(V),\,2d_{H_2}(V),\,4d_{K_4}(V)\big),\\
&k(V)=\max\big(n,\,2d_{H_0}(V)-4,\,2d_{H_1}(V)-4,\,2d_{H_2}(V)-4,\,4d_{K_4}(V)-12\big).
\end{align*}
In particular, $X(V):= S^V\wedge H\uF$ is an $n$-slice if and only if $\lceil n/2\rceil\le d_{H_i}(V)\le\lfloor n/2\rfloor+2$ for $i=0,1,2$.
\end{cora}
See Corollary~\ref{pk4range}.

Note that $2d_{K_4}(V)=d_{H_0}(V)+d_{H_1}(V)+d_{H_2}(V)-n$. Thus the conditions  $4d_{K_4}(V) \ge n$ and  $4d_{K_4}(V)-12\le n$  automatically follow from the conditions on $d_{H_i}(V)$ for $i=0,1,2$.

\subsection{Cofibres of Thom classes}
Different $X(V)$ are connected by Thom classes. For the sign representation $\alpha_i$, the class $u_{\alpha_i}$ gives a map
$X(V+1-\alpha_i)\to X(V)$ whose cofibre
$C(V,\alpha_i)$ is supported on the subgroups $H_j$ and $K_4$ for $j\neq i$. Every comparison between slice towers of the spheres $X(V)$ passes through
these cofibres, so we determine their slices completely.

The cofibre is invisible to $e$ and $H_i$, and it is built from three simpler pieces: a $H_j$-part for $j\neq i$, each pulled back
from a $C_2$-spectrum, glued along a spectrum seen only by $K_4$. Our next result shows that the slice tower respects this
decomposition.

\begin{thmb}\label{tB}
The cofibre \(C(V,\alpha_i)\) admits only even-dimensional slices, which are explicitly described in Theorem~\ref{tCslices}. 
\end{thmb}


\subsection{The slices of $X(\infty\alpha_i)$}
The spectrum $X(\infty\alpha_i)=S^{\infty\alpha_i}\wedge H\uF$
is the colimit of the spectra $X(n\alpha_i)$ along multiplication
by the Euler class $a_{\alpha_i}$ of the representation $\alpha_i$, and it controls the slices of $X(n\alpha_i)$ below level $n$. Its slice tower can be described completely. We use $[2]$ to denote the set $\{0,1,2\}$.
\begin{thmc}\label{tC}
The slices of $X(\infty \alpha_i)$ are given by
 \[
P^{4m}_{4m}X(\infty \alpha_i)\cong
\begin{cases}
 H\ug_{(i)} & m=0,\\
S^{ \rho}\wedge  H\ug_{(i)}^{*} & m=1,\\
S^{m \rho}\wedge  H(\phi^*_{H_j}\uF^{*}\oplus\phi^*_{H_k}\uF^{*}\oplus\ug^{\oplus(2m-3)}) & m\ge2.
\end{cases}
\]
and 
\[
P^{4m+2}_{4m+2}X(\infty\alpha_i)\simeq
\begin{cases}
\Sigma H\ug_{(i)} & m=0,\\
S^{m\rho}\wedge\Sigma H(\phi^{\ast}_{H_j}\uf_{C_2}\oplus\phi^{\ast}_{H_k}\uf_{C_2}) & m\ge1,
\end{cases}
\]
where $i,j,k$ are distinct elements of the set $[2]$.
\end{thmc}
See Theorem~\ref{thm:Xinfslices}.

The proof combines $\rho$-periodicity, Slone's classification of $2$-slices for $K_4$, and the $RO(K_4)$-graded Bredon
cohomology of certain universal space. Under Brown--Comenetz duality, the restriction maps that detect the slices become transfers, and
these can be computed by Frobenius reciprocity.

\medskip

\subsection{Defects and spherical slices}~
We return to the windows of Corollary~A.1. What happens when $V$ ($n=\dim V$) lies outside the window? Corollary~A.1. suggests a simple strategy. Each subgroup $H$ of order two sees a window $[\lceil n/2\rceil,\lfloor n/2\rfloor+2]$ for its fixed-point dimension $d_H(V):= \dim(V^H)$. If $d_H(V)$ falls outside, the most economical correction is to push it back to the nearest admissible value.
Let $\delta_H(V)$ be the nearest point of $[\lceil n/2\rceil,\lfloor n/2\rfloor+2]$ to $d_H(V)$. The \emph{defect vector} of $V$ are
\[
\Delta(V)=(\Delta_0(V),\Delta_1(V),\Delta_2(V)),\qquad \Delta_{H_i}(V)=\delta_{H_i}(V)-d_{H_i}(V)\quad(i=0,1,2).
\]
If $\sum_{\Delta(V)}:=\Delta_0(V)+\Delta_1(V)+\Delta_2(V)$ is even, then define
\[
x_i(V)=\tfrac12(\Delta_j(V)+\Delta_k(V)-\Delta_i(V))
\]
 for $\{j,k\}=[2]\setminus\{i\}$.
 
The defect vector measures how far each fixed-point dimension lies outside the window of Corollary~\ref{cA}. In particular, $\Delta(V)=0$ exactly when $X(V)$ is an $n$-slice. 

Let $V$ be a $G$-representation and $\uM$ a Mackey functor. We say that $S^V\wedge H\uM$ has a \emph{spherical slice} if $P^{n}_{n}(S^V\wedge H\uM)\simeq S^U\wedge H\uM$ for some
$U\in RO(G)$ with $\dim U=n$.
\begin{thmd}\label{tD}
Suppose $\sum_{\Delta (V)}$ is even.
\begin{enumerate}
\item If $X(V)\le n$ and $x_0(V),x_1(V),x_2(V)\ge0$, then $X(V)$ has a spherical slice.
\item If $X(V)\ge n$ and $x_0(V),x_1(V),x_2(V)\le0$, the same conclusion holds.
\end{enumerate}
\end{thmd}
See Theorem~\ref{main}.
 For the other case, $\sum_{\Delta(V)}$ is odd.
\begin{thme}
Suppose $\sum_{\Delta(V)}$ is odd.
\begin{enumerate}
\item  Let $X(V)\le n$ and
$x_0(V),x_1(V),x_2(V)\ge0$. Assume $\Delta_i(V)$ is odd for some $i\in [2]$. Then $P^n_nX(V)$ is not spherical, and there is a cofibre sequence
\[
X(\widehat{V'})\xrightarrow{\ v_i\ }P^n_nX(V)\to P_n(\Sigma^{\widehat{V'}-\tau_0}C(v_i)),
\]
in which $X(\widehat{V'})$ is an $n$-slice for some $\widehat{V'}.$ 
\item If $X(V)\ge n$ and $x_0(V),x_1(V),x_2(V)\le0$. Assume $\Delta_i(V)$ is odd.  Then $P^n_nX(V)$ is not spherical, and there is a cofibre sequence
\[
P^n(\Sigma^{\widehat{V'}-1}C(v_i))\to  P^n_nX(V) \xrightarrow{\ v_i\ } X(\widehat{V'}),
\]
in which $X(\widehat{V'})$ is an $n$-slice for some $\widehat{V'}.$
\end{enumerate}
\end{thme}
See Theorem~\ref{mainodd}.

\begin{thmf}
    (a) Let $n\ge 5$. Then $X(n) \ge n$, $\Delta(n)= (-k, -k, -k)$ and $x_i=-\frac{k}{2}$, where $i=0,1,2$ and $k = \lceil n/2\rceil-2.$ 
    \begin{enumerate}
        \item If $n=0,3 \pmod 4$,  then $P^n_nX(n)$ is spherical and $P^n_n X(n)\simeq X(n-\frac{k}{2}(3-\bar{\rho})).$
        \item If $n=1,2 \pmod 4,$ then $P^n_n X(n)$ is not spherical and there is cofibre sequence
        \[
         P^n(\Sigma^{\widehat{V'}-1}C(v_0)) \to P^n_n X(n) \xrightarrow{\ v_0\ }X\bigl(\widehat{V'}\bigr) 
        \]
        where $\widehat{V'}=(n-\frac{\lceil n/2\rceil-1}{2}-k)+(k-\frac{\lceil n/2\rceil-1}{2})\,\alpha_0+\frac{\lceil n/2\rceil-1}{2}\,\alpha_1+\frac{\lceil n/2\rceil-1}{2}\,\alpha_2$ and $X(\widehat{V'})$ is an $n$-slice.
    \end{enumerate}

    Moreover, the $n$-slice is
    \[
    P^n_nX(n)\simeq S^{y\rho+(n-4y)} \wedge H\uf,
    \]
    where $y=\frac{\lceil n/2\rceil-1}{2}.$

    (b) For $0\leq n\le 4$, $X(n)$ is an $n$-slice.
\end{thmf}
See Theorem~\ref{GYThm}.

\begin{mysubsection}{Organization}~We organize the remainder of the paper as follows. In Section~2, we present the preliminaries needed for the subsequent sections. In Section~3, we describe the slice interval of the Mackey functor $\uM$ for a finite group $G$ (Theorem~\ref{trange}). In particular, for $G=K_4$ and the constant Mackey functor $\uF$, we obtain Corollary~\ref{pk4range}. In Section~4, we study the cofibres of certain maps. Specifically, we describe the cofibres of the maps induced by the Thom class $u_{\alpha_i}$ and a certain twist class $v_i$. In Section~5, we determine the fibre and its slices of the map $S^{n\alpha_i}\to S^{\infty\alpha_i}$ in Theorem~\ref{thm:Xinfslices} and consequently, compute all the slices of $S^{\infty\alpha_i}$. Finally, in Section~6, we introduce the defect vector and the notion of spherical slices. We characterize spherical slices for different cases of defect vectors in Theorems~\ref{main} and~\ref{mainodd}. At the end of this section, we provide an alternative proof of the result of~\cite{GY20}, stated in Proposition~\ref{GYThm}. In the appendix, we describe all the Mackey functors used throughout this article. 
\end{mysubsection}

\paragraph{\textbf{Acknowledgement.}}The research of Surojit Ghosh was supported by ANRF, India, with grant number  ARG/2025/000434/MS. Ankit Kumar's research is supported by the University Grants Commission, India.

\begin{subsection}{Notation}
We fix the following notation throughout the paper.
\begin{enumerate}
\item We use $G$ to denotes a finite group, and $K_4$ to denote the Klein-four group.
\item Let \(H\subseteq G\). For any \(V\in RO(G)\), define $d_H(V):=\dim(V^H)$. When the representation \(V\) is clear from the context, we write \(d_H\) for \(d_H(V)\).
\item Throughout, we use $[2]$ to denote the set \(\{0,1,2\}\).
\item Let $G$ be a finite group and $\uM$ a Mackey functor. For $V\in RO(G)$, we write
\[
X(V,\uM):=S^V\wedge H\uM.
\]
When $G=K_4$ and $\uM=\uF$ is the constant Mackey functor, we abbreviate $X(V,\uF)$ to $X(V)$ for $V\in RO(K_4)$.
\item We write $\rho_{C_2}:=1+\sigma$ for the regular representation of $C_2$, where $\sigma$ denotes the sign representation. We use $\rho$ to denote the regular representation of $K_4$.
\item We write $\sum_ {\Delta_{(V)}}$ for $\Delta_0(V)+\Delta_1(V)+\Delta_2(V)$.
\end{enumerate}
\end{subsection}

\section{Background}
In this article, we use $K_4$ (generated by $a,b$) to denote the Klein four-group. Throughout, $H_0=\langle a\rangle,  H_1=\langle{ab}\rangle,  H_2=\langle b\rangle,$ are the three subgroups of order two. Let $\aO,\aI,\aII$ be the three sign representations, normalised by $\ker\aO=H_0, \ker\aI=H_1,\ker\aII=H_2,$
and the regular representation  $\rho=1+\aO+\aI+\aII$. Restricting to the order-two subgroups (and writing $\sigma$ for the sign representation of $C_2$):
\begin{myeq}\label{erest}
\begin{aligned}
(a+b\aO+c\aI+d\aII)|_{H_0}&=(a+b)+(c+d)\sigma,\\
(a+b\aO+c\aI+d\aII)|_{H_1}&=(a+c)+(b+d)\sigma,\\
(a+b\aO+c\aI+d\aII)|_{H_2}&=(a+d)+(b+c)\sigma.
\end{aligned}
\end{myeq}

\begin{defn}
    The \emph{real representation ring} of a group $G$ is the free abelian group generated by all irreducible representations. In particular, $RO(K_4):=\mathbb Z\{1,\aO,\alpha_1,\alpha_2\}.$  
\end{defn}
\begin{mysubsection}{Mackey Functors}
~  Ordinary cohomology theories are defined for Abelian groups, and these are represented by spectra with homotopy concentrated in degree $0$. In the equivariant world, the analogous role is played by \emph{Mackey functors}. In this section we briefly recall their definition, and relate them to equivariant cohomology (see \cite{May96} for details). 

The \emph{Burnside category} $\Burn_G$ is the category with objects finite $G$-sets
 and each morphism set $\Mor_{\Burn_G}(S,T)$ is the group completion of the hom-set of spans between $S$ and $T$ in the category of finite $G$-sets.  $\Burn_G$ is self-dual; that is, the duality map $D\colon  \Burn_G \to \Burn_G^{\op}$ that is identity on objects and switches the legs of the spans, is an isomorphism of categories.
 
 \begin{defn}
A product-preserving functor $\uM: \Burn_G^{op} \to \Ab$ from the Burnside category into abelian groups is called a \emph{Mackey functor}.
\end{defn}

In this paper we restrict our attention to $G=K_4$. 
Explicitly, a $G$-Mackey functor\footnote{This is a simplification in the case $G$ is Abelian. Otherwise the double coset formula (4) has a slightly more complicated expression.} $\uM$ is a collection of commutative $W_G(H)$-groups $\uM(G/H)$ one for each subgroup $H \le G$, each accompanied by \emph{transfer} $\tr^H_K\colon  \uM(G/K) \to \uM(G/H)$ and \emph{restriction} $\res^H_K \colon\uM(G/H) \to \uM(G/K)$ for $K\le H\le G$ such that
\begin{enumerate}
\item $\tr^H_J = \tr^H_J\tr^K_J$ and $\res^H_J = \res^K_J \res^H_K$ for all $J \le K \le H.$
\item  $\tr^H_K(\gamma.x)= \tr^H_K(x)$ for all $x \in \uM(G/K)$ and $\gamma \in W_H(K).$
\item $\gamma. \res^H_K(x) =\res^H_K(x)$ for all $x \in \uM(G/H)$ and $\gamma \in W_H(K).$
\item $\res^H_K\tr ^J_K(x)= \sum_{\gamma \in W_H(K)} \gamma.\tr^{K}_{J\cap K}(x)$
for all subgroups $J,H \leq K.$
\end{enumerate}

Mackey functors play a fundamental role in equivariant stable homotopy theory, particularly in describing the equivariant homotopy groups of \(G\)-spectra. For $G=K_4$, we represent a general Mackey functor $\uM$ by its Lewis diagram~\cite{L88}:

\[
\xymatrix@C=3.1em@R=2.5em{
 & \uM(K_4/K_4)
   \ar@/_0.5pc/[ld]_{\res}
   \ar@/^0.5pc/[d]^{\res}
   \ar@/^0.5pc/[rd]^{\res}
 & \\
\uM(K_4/H_0)
   \ar@/_0.5pc/[ru]_{\tr}
   \ar@/_0.5pc/[rd]_{\res}
 & \uM(K_4/H_1)
   \ar@/^0.5pc/[u]^{\tr}
   \ar@/^0.5pc/[d]^{\res}
 & \uM(K_4/H_2)
   \ar@/^0.5pc/[lu]^{\tr}
   \ar@/^0.5pc/[ld]^{\res}
 \\
 & \uM(K_4/e).
   \ar@/_0.5pc/[lu]_{\tr}
   \ar@/^0.5pc/[u]^{\tr}
   \ar@/^0.5pc/[ru]^{\tr}
 &
}
\]
Here, ${H_0}$, ${H_1}$, and ${H_2}$ denote the left, diagonal, and right cyclic subgroups of order two, respectively, while $e$ denotes the trivial subgroup. The downward arrows represent restriction maps, and the upward arrows represent transfer maps between the corresponding values of the Mackey functor on the associated orbits.

\begin{example} A basic example is the \emph{constant Mackey functor} \(\underline{C}\), associated to an abelian group \(C\), defined by \(\underline{C}(G/H)=C\) for every subgroup \(H\leq G\). The restriction maps are identities, while the transfer maps \(\operatorname{tr}_K^H\), for \(K\leq H\), are given by multiplication by the index \([H:K]\). 
\end{example}
A selection of Mackey functors relevant to the present article is provided in the Appendix.
\end{mysubsection}

Let \(\Sp^G\) denote the category of orthogonal \(G\)-spectra \cite{MM02}. It admits a symmetric monoidal model category structure \((\Sp^G,\wedge,\mathbb{S})\), as described in \cite{HHR16}, and we use it as our model for \(G\)-equivariant stable homotopy theory.

Equivariant homotopy groups are naturally graded by \(RO(G)\), the real representation ring of \(G\). For \(V=W-U\in RO(G)\), we define
\[
\pi_V^G(X):=[S^W,S^{U}\wedge X]^G,
\]
where the right-hand side denotes the group of homotopy classes of maps in the \(G\)-equivariant stable homotopy category. When the group $G$ is clear from the context, we omit it from the notation and write $\pi_V(X)$ instead. Thus, every orthogonal \(G\)-spectrum \(X\) has \(RO(G)\)-graded homotopy groups, denoted by \(\pi_\star (X)\). These homotopy groups assemble into a Mackey functor \(\underline{\pi}_V(X)\), whose value at \(G/H\) is
\[
\underline{\pi}_V(X)(G/H)
:=\pi_V(\Sigma^{\infty}_{G}G/H_{+}\wedge X).
\]

For every \(G\)-Mackey functor \(\underline{M}\), there exists an equivariant Eilenberg--Mac Lane spectrum \(H\underline{M}\) \cite{GM95}, characterized by
\[
\upi_k(H\underline{M})\cong
\begin{cases}
\underline{M} & k=0,\\
0 & k\ne 0,
\end{cases}
\]
where the vanishing condition is understood in the appropriate equivariant homotopy groups. The spectrum \(H\underline{M}\) represents a reduced cohomology theory given by
\[
\widetilde{H}_G^V(X;\underline{M})
\cong [\Sigma^\infty X,\Sigma^V H\underline{M}]^G.
\]
This theory is called \(RO(G)\)-graded Bredon cohomology with coefficients in the Mackey functor \(\underline{M}\).

\begin{mysubsection}{Some special classes in Bredon cohomology} 
~The $RO(G)$-graded cohomology ring $\widetilde{H}^{\bigstar}_{G}(S^0; \uF)$ contains distinguished classes associated to representations \cite{HHR16}. For any real $G$-representation $V$ with no nontrivial fixed point, there is a canonical \emph{Euler class}
    \[
    a_V \in \widetilde{H}_G^V(S^0; \uF),
    \]
    represented by the inclusion $S^0 \hookrightarrow S^V$.
If $V$ is oriented representation, there exists a \emph{Thom class}
    \[
    u_V \in \widetilde{H}_G^{V - \dim V}(S^0; \uF),
    \]
    arising from the cofibre sequence $S(V)_{+} \to D(V)_{+} \to S^V$.

These classes satisfy natural multiplicative identities:
\[
a_V \cdot a_W = a_{V \oplus W}, \quad u_V \cdot u_W = u_{V \oplus W}.
\]

\begin{prop}[\cite{BG21}, \cite{GK26}]{\label{coh2}}
    \begin{enumerate}
        \item  The Mackey functor $\underline{H}^{\bigstar}_{K_4}(S^0; \uF)$ is given by
\[
\vcenter{\xymatrix@R=1.2cm{
\mathbb{M}_2(\sigma):=\F[a_\sigma,u_\sigma]\ \oplus\ \displaystyle\bigoplus_{j,k\ge1}\F\Bigl\langle\Sigma^{-1}\frac{1}{a_\sigma^ju_\sigma^k}\Bigr\rangle
 \ar@/_1.2pc/[d]_{\res}\\
\F[u^{\pm1}]\ar@/_1.2pc/[u]_{\mathrm{tr}}
}}
\]
 with $\res(a_\sigma)=0$, $\res(u_\sigma)=u$, $\tr(u^k)=0$ for $k\ge -1$ and $\tr(u^{-m})=\Sigma^{-1}\frac{1}{a_\sigma u_\sigma^{\,m-1}}$ for $m\ge2$. Here $\sigma$ denotes the sign representation of the group $C_2$.
        
\item Consider the $K_4$-universal space $E\mathcal F$ for the family $\mathcal F=\{e,H_0\}$. Then the Mackey functor
$\underline{H}^{\bigstar}_{K_4}(E\mathcal F_+;\uF)$ is given by
\[
\xymatrix@C=0.9cm@R=1.4cm{
 & \mathbb{M}_2(\aII)\otimes\F[a_{\aO},u^{\pm}_{\aO},u^{\pm}_{\alpha_1-\alpha_2}]
   \ar@/_0.5pc/[ld]_{\res_{H_0}}
   \ar@/^0.5pc/[d]^{\res_{H_1}}
   \ar@/^0.5pc/[rd]^{\res_{H_2}} & \\
\mathbb{M}_2(\sigma)
 \ar@/_0.5pc/[ru]_{\tr^{K_4}}
   \ar@/_0.5pc/[rd]_{\res_e}
 & \F[a_\sigma,u_\sigma^{\pm}]
   \ar@/^0.5pc/[u]^{\tr^{K_4}}
   \ar@/^0.5pc/[d]^{\res_e}
 & \F[a_\sigma,u_\sigma^{\pm}]
    \ar@/^0.5pc/[lu]^{\tr^{K_4}}
   \ar@/^0.5pc/[ld]^{\res_{e}}\\
 & \F[u^{\pm}_{\aO},u^{\pm}_{\aI},u^{\pm}_{\aII}].
    \ar@/_0.5pc/[lu]_{\tr^{H_0}}
   \ar@/^0.5pc/[u]^{\tr^{H_1}}
   \ar@/^0.5pc/[ru]^{\tr^{H_2}}&
}
\]  
The restrictions from $K_4$ are determined by
\begin{align*}
\res_{H_0}&:\ a_{\aO}\mapsto0,\quad u_{\aO}\mapsto1,\quad a_{\aII}\mapsto a_\sigma,\quad u_{\aII}\mapsto u_\sigma,\quad u_{\alpha_1-\alpha_2}\mapsto1,\\
\res_{H_1}&:\ a_{\aO}\mapsto a_\sigma,\quad u_{\aO}\mapsto u_\sigma,\quad a_{\aII}\mapsto a_\sigma,\quad u_{\aII}\mapsto u_\sigma,\quad u_{\alpha_1-\alpha_2}\mapsto u_\sigma^{-1},\\
\res_{H_2}&:\ a_{\aO}\mapsto a_\sigma,\quad u_{\aO}\mapsto u_\sigma,\quad a_{\aII}\mapsto0,\quad u_{\aII}\mapsto1,\quad u_{\alpha_1-\alpha_2}\mapsto u_\sigma .
\end{align*}
    \end{enumerate}
\end{prop}
\begin{proof}(1) The restrictions of all these classes are determined by degree considerations and the ring structure. On the other hand, the transfer \(\tr(u^k)\) vanishes for \(k=-1\) by degree considerations and for \(k\geq 0\) by the Frobenius reciprocity identity.

To determine \(\tr(u^{-m})\), consider the following exact sequence induced by the cofibre sequence \(C_2/e_+\to S^0\to S^\sigma\):
\[
0 \longrightarrow
\pi_{1}(\Sigma^{m-m\sigma}H\uF)^{C_2}
\xrightarrow{f}
\pi_{1}(\Sigma^{m-(m-1)\sigma}H\uF)^{C_2}
\xrightarrow{\partial}
\pi_{0}(\Sigma^{m-m\sigma}H\uF)^{e}
\xrightarrow{\tr}
\pi_{0}(\Sigma^{m-m\sigma}H\uF)^{C_2}
\longrightarrow 0.
\]
When \(m=2\), the domain of \(\partial\) is zero, so exactness implies that \(\tr\) is an isomorphism. For \(m\geq 3\), the map \(f\) is an isomorphism. Consequently, \(\partial=0\), and exactness again implies that \(\tr\) is an isomorphism. 

    (2) The description of the \(RO(K_4)\)-graded Bredon cohomology of \(E\cF\) follows from \cite{GK26}, where it is done using the homotopy fixed point spectral sequence. For the computation of $\res_{H_0}$ and $\res_{H_\ell}$ for $\ell=1,2$, consider the following diagrams 
\[
\begin{array}{c@{\quad\quad}c}
\xymatrix{
K_4/K_4: & \widetilde{H}^{0}_{K_4}(E\cF_{+})
\ar[r]^{.a_{\alpha_0}} \ar[d]_{\res_{H_0}}
& \widetilde{H}^{\alpha_0}_{K_4}(E\cF_{+})
\ar[d]^{\res_{H_0}} \\
K_4/H_0: & \widetilde{H}^{0}_{H_0}(S^0)
\ar[r]_0
& \widetilde{H}^{1}_{H_0}(S^0)\cong 0
}
&
\xymatrix{
K_4/K_4: & \widetilde{H}^{0}_{K_4}({E\cF}_{+}) \ar[r]^{.a_{\alpha_0}}\ar[d]_{\res_{H_\ell}}
& \widetilde{H}^{\alpha_0}_{K_4}({E\cF}_{+}) \ar[d]^{\res_{H_\ell}} \\
K_4/H_\ell: & \widetilde{H}^{0}_{H_\ell}({EC_2}_{+}) \ar[r]_{.a_{\sigma}} & \widetilde{H}^{\sigma}_{H_\ell}({EC_2}_{+})
}
\end{array}
\]
Using the commutativity of the above diagram, we can determine the restriction of \(a_{\alpha_0}\). Replacing multiplication by \(a_{\alpha_0}\) with multiplication by \(a_{\alpha_\ell}\) yields the restrictions of the other \(a\)-classes. The restrictions of all remaining classes which are mention are determined by the fact that their restrictions to the trivial subgroup is $1$.
    \end{proof}

\begin{lemma}\label{kerequiv}
    We have the following equivalences in $\Sp^{K_4}:$
    \begin{enumerate}
        \item $\Sigma^{\alpha_i-1}H\uf \simeq H\ker(\uF \twoheadrightarrow \phi^*_{H_i}\uF)$.
        \item For $i\neq j$, $S^{\alpha_i -\alpha_j}\wedge H\ker(\uF \twoheadrightarrow \phi^*_{H_j}\uF) \simeq H\ker(\uF \twoheadrightarrow \phi^*_{H_i}\uF)$.
    \end{enumerate}
\end{lemma}
\begin{proof}
    The first equivalence follows from the cofibre sequence
    \[
    {K_4/H_i}_+ \to S^0\to S^{\alpha_i}
    \]
    and smashing with $\Sigma^{-1}H\uf.$

    Whereas the second one is the consequences of the first.
\end{proof}
\end{mysubsection}

\begin{mysubsection}{$\uZ$-modules and Inflation  functors}
~Let $\BZ_G$ be the additive category whose objects are finite $G$-sets and whose morphisms are
\[
\mathrm{Mor}_{\BZ_G}(S,T):=\mathrm{Hom}_{\mathbb Z[G]}\big(\mathbb Z[S],\mathbb Z[T]\big),
\]
i.e.\ the full subcategory of permutation modules. 

Note that 
\begin{myeq}\label{eq:homfix}
\mathrm{Hom}_{\mathbb Z[G]}\big(\mathbb Z[G/H],\mathbb Z[T]\big)\xrightarrow{\ \cong\ }\mathbb Z[T]^H,\qquad f\mapsto f(eH).
\end{myeq}

By \cite[Prop.~2.15]{Zen17} (a form of Yoshida's theorem), $\uZ$-modules for $G$ are the same as cohomological Mackey functors. 
We write $\ZMod_G$ for this category.

The quotient $q: G \to G/N$ induces a functor
\[
\phi\colon(\BZ_{G/N})^{\rm op}\to(\BZ_G)^{\rm op} \qquad T\mapsto q^*T,
\]
which regards a $G/N$-set as a $G$-set.
 This gives a functor $\phi^*: \ZMod_G  \to \ZMod_{G/N}$ given by precomposition, that is, it sends $\uM \mapsto \uM \circ \phi.$
\begin{defn}
For $\uM\in\ZMod_{G/N}$, define
\[
\Phi^*_N\uM:=\Lan_\phi(\uM)\in\ZMod_G ,
\]
the (additive) left Kan extension of $\uM$ along $\phi$. Concretely, for a finite $G$-set $X$,
\begin{myeq}\label{eq:coend}
\Phi_N^*\uM(X)= \operatorname*{colim}_{(X\to\phi T)}\uM(T),
\end{myeq}
where the colimit runs over the comma category of $\BZ_G$-morphisms $X\to\phi T$.
\end{defn}
 By the definition of the left Kan extension, the functor $\Phi_N^*$ is left adjoint to $\phi^\ast.$
\begin{lemma}\label{lemkey}
For $H\le G$ and $T\in\BZ_{G/N}$,  there is an isomorphism, natural in $T$,
\[
\mathrm{Hom}_{\BZ_G}\big(G/H,\phi T\big)\;\cong\;\mathrm{Hom}_{\BZ_{G/N}}\big((G/N)/\Hb,\,T\big)
\]
where $\Hb$ is the image of $H$ in $G/N$.
\end{lemma}
\begin{proof}
By \eqref{eq:homfix} the left side is $\mathbb Z[T]^H$. Since $N$ acts trivially on $T$, $H$ acts through its image $\Hb$, so $\mathbb Z[T]^H=\mathbb Z[T]^{\Hb}$, which is the right side by \eqref{eq:homfix} for $G/N$. Both identifications are natural in $\mathbb Z[G/N]$-linear maps $\mathbb Z[T]\to\mathbb Z[T']$. Under them, $f\leftrightarrow x=f(eH)$ and $\bar f\leftrightarrow x$, and $\phi(\bar f)\circ p_H$ sends $eH\mapsto eHN\mapsto x$.
\end{proof}

\begin{cor}\label{tformula}
For $\uM\in\ZMod_{G/N}$ and $H\le G$,
\[
\Phi^*_N\uM(G/H)\;\cong\;\uM\big((G/N)/(HN/N)\big)=
\begin{cases}
\uM\big((G/N)/(H/N)\big) & \text{if } N\subseteq H,\\[2pt]
\uM\big((G/N)/(HN/N)\big) & \text{if } N\not\subseteq H .
\end{cases}
\]
In particular, $\Phi_N^*\uM(G/H)=\uM\big((G/N)/e\big)$ whenever $H\subseteq N$.
\end{cor}
\begin{proof}
By Lemma~\ref{lemkey}, the functor $T\mapsto\mathrm{Hom}_{\BZ_G}(G/H,\phi T)$ is represented by $(G/N)/\Hb$, so by the co-Yoneda lemma the coend \eqref{eq:coend} is $\uM((G/N)/\Hb)$. 
\end{proof}


For the above functor $\Phi_N^{*}$, we have the following example, which we will use in the subsequent sections 
\[
\ug_{(0)}=\Phi^*_{H_0}\ug_{C_2},\qquad \ug_{(1)}=\Phi^*_{H_1}\ug_{C_2},\qquad \ug_{(2)}=\Phi^*_{H_2}\ug_{C_2} .
\]

Their corresponding Lewis diagrams are given in Example~\ref{Phimac} in Appendix.

\begin{cor}\label{cprops}
\begin{enumerate}
\item The functor \(\Phi_N^*\) is exact, as follows from its pointwise description, and is left adjoint to the functor
\[
q_*\colon \mathsf{Mack}_G\longrightarrow \mathsf{Mack}_{G/N},
\]
defined by
\[
(q_*\underline M)((G/N)/(H/N))=\underline M(G/H).
\]

\item The restriction \(\operatorname{res}_N^G\Phi_N^*\underline M\) is the constant \(\underline{\mathbb Z}\)-module associated to \(\underline M((G/N)/e)\), with identity restriction maps and transfer maps given by multiplication by the subgroup index.

\item If \(N=e\), then \(\Phi_e^*=\mathrm{id}\). If \(N=G\), then \(\Phi_G^*A=\underline A\) is the constant Mackey functor. In particular,
\[
\Phi_N^*\underline{\mathbb Z}=\underline{\mathbb Z}
\text{  and  }
\Phi_N^*\uF=\uF.
\]

\item At the spectrum level, \(\Phi_N^*\) extends to a functor on \(H\underline{\mathbb Z}\)-modules, given by
\[
X\longmapsto
H_G\underline{\mathbb Z}
\wedge_{q^*H_{G/N}\underline{\mathbb Z}}q^*X.
\]
It satisfies
\[
\Phi_N^*(S^V\wedge X)
\simeq S^{q^*V}\wedge\Phi_N^*X.
\]
Moreover, when \(\mathsf{Mack}_{G/N}\) has finite projective dimension (for example, when \(G/N=C_2\)), we have
\[
\Phi_N^*H\underline M
\simeq H\Phi_N^*\underline M
\text{ and }
\underline\pi_k\Phi_N^*X
\cong\Phi_N^*\underline\pi_kX.
\]
\end{enumerate}
\end{cor}
\begin{proof}
(1) Kan extension adjunction; exactness from Theorem~\ref{tformula}. 

 (2) For $H'\le H\le N$, $\Hb=\overline{H'}=e$ and the transfer factor is $[H:H']$. 
 
 (3) Formal from the formula. 
 
 (4) As in \cite[\S 3]{GS25}: check on projectives $\uparrow\uZ$ and use exactness.
\end{proof}
\end{mysubsection}

\begin{mysubsection}{Geometric inflation $\phi^*_N$ and inflation by zero}
The geometric inflation $\phi^*_N:=\widetilde{E\mathcal F[N]}\wedge q^*(-)$ is right adjoint to the geometric fixed points $\Phi^N$; it is fully faithful onto $\widetilde{E\mathcal F[N]}$-local spectra, and $\Phi^Nq^*\simeq\mathrm{id}$. On Mackey functors its shadow is inflation by zero:
\[
\phi^*_N\underline M(G/H)=\begin{cases}\underline M((G/N)/(H/N))&N\subseteq H,\\0&N\not\subseteq H,\end{cases}\qquad \phi^*_NH\underline M\simeq H\phi^*_N\underline M .
\]

For $i=0,1,2$, let $q\colon K_4\to K_4/H_i\cong C_2$, so that $q^*\sigma=\alpha_i$ by definition of $\alpha_i$. On $H\uZ$-modules,
\[
\Phi^*_{H_i}Z:=H\uZ\wedge_{q^*H_{C_2}\uZ}q^*Z ,
\]
using the ring map $\eta\colon q^*H_{C_2}\uZ\to H\uZ$ adjoint to the equivalence $H_{C_2}\uZ\simeq(H\uZ)^{H_i}$.

\begin{lemma}\label{lem:Phibasic} For $G=K_4$ and $i=0,1,2$.
\begin{enumerate}
\item $\Phi^*_{H_i}H_{C_2}\uZ\simeq H\uZ$.
\item $\Phi^*_{H_i}$ is exact; therefore it preserves cofibre sequences and all colimits.
\item $\Phi^*_{H_i}H_{C_2}\uF\simeq H\uF$.
\item For every $V\in RO(C_2)$ and every $H_{C_2}\uZ$-module $Z$, $\Phi^*_{H_i}(S^V\wedge Z)\simeq S^{q^*V}\wedge\Phi^*_{H_i}Z$.
\end{enumerate}
\end{lemma}
\begin{proof}
(a) Base change of the unit module: $H\uZ\wedge_{q^*H\uZ}q^*H\uZ\simeq H\uZ$.

(b) $q^*$ is a left adjoint between stable categories, and so is base change $H\uZ\wedge_{q^*H\uZ}(-)$.

(c) $H_{C_2}\uF$ is the cofibre of $2\colon H_{C_2}\uZ\to H_{C_2}\uZ$ in $H_{C_2}\uZ$-modules, from $0\to\uZ\xrightarrow{2}\uZ\to\uF\to0$. By (a) and (b), $\Phi^*_{H_i}H_{C_2}\uF$ is the cofibre of $\Phi^*_{H_i}(2)=2\colon H\uZ\to H\uZ$, which is $H\uF$ by $0\to\uZ\xrightarrow2\uZ\to\uF\to0$ over $K_4$.

(d) $q^*$ is strong symmetric monoidal with $q^*S^V=S^{q^*V}$ for actual $V$, hence for virtual $V$ (the spheres are invertible). Therefore $q^*(S^V\wedge Z)\simeq S^{q^*V}\wedge q^*Z$ as $q^*H\uZ$-modules, and
\[
H\uZ\wedge_{q^*H\uZ}\big(S^{q^*V}\wedge q^*Z\big)\simeq S^{q^*V}\wedge\big(H\uZ\wedge_{q^*H\uZ}q^*Z\big).\qedhere
\]
\end{proof}

\begin{prop}\label{pXninfl}
For every $n\in\mathbb Z$ and $K_4$-Mackey functor $\underline M\in\{\uZ,\uF\}$, we have
\begin{myeq}{\label{PhiL}}
X_n(\underline M)=S^{n\alpha_i}\wedge H\underline M\ \simeq\ \Phi^*_{H_i}\big(S^{n\sigma}\wedge H_{C_2}\underline M\big)\text{ for } i=0,1,2.
\end{myeq}
More generally $S^{a+b\alpha_i}\wedge H\underline M\simeq\Phi^*_{H_i}\big(S^{a+b\sigma}\wedge H_{C_2}\underline M\big)$ for all $a,b\in\mathbb Z$. In particular:
\begin{enumerate}
\item $X_n(\underline M)^{H_i}\simeq S^{n\sigma}\wedge H_{C_2}\underline M$.
\item $\res^{K_4}_{H_j}X_n(\underline M)\simeq S^{n\sigma}\wedge H_{C_2}\underline M$ and $\res^{K_4}_{H_i}X_n(\underline M)\simeq S^n\wedge H_{C_2}\underline M$ for $j\neq i$.
\end{enumerate}
\end{prop}
\begin{proof}
By Lemma~\ref{lem:Phibasic}(d) with $V=a+b\sigma$, then (a) resp.\ (c),
\begin{align*}{\label{PhiL}}
\Phi^*_{H_i}\big(S^{a+b\sigma}\wedge H_{C_2}\underline M\big)\simeq S^{q^*(a+b\sigma)}\wedge\Phi^*_{H_i}H_{C_2}\underline M\simeq S^{a+b\alpha_i}\wedge H\underline M .    
\end{align*}

$(1)$ Then follows from $(\Phi^*_{H_i}Z)^{H_i}\simeq Z$.

$(2)$ Since $\alpha_i|_{H_j}=\sigma$ and  $\alpha_i|_{H_i}$ is trivial.
\end{proof}

\end{mysubsection}

\section{The slice interval}
The regular slice filtration is a fundamental tool in equivariant stable homotopy theory, providing an analogue of the Postnikov filtration that organizes equivariant spectra according to their representation-theoretic dimensions. It decomposes a \(G\)-spectrum into layers, called \emph{slices}, which capture equivariant information at different levels and are particularly useful in computing equivariant homotopy groups and studying equivariant phenomena.
\begin{mysubsection}{The slice tower}
~Let $\tau_{\ge n}\subseteq\Sp^G$ be the localizing subcategory generated by the
regular slice cells of dimension $\ge n$ \cite[\S4]{HHR16}. Associated to $\tau_{\ge n}$
there is a natural localization functor $P^{n-1}$, whose acyclics are exactly the objects
of $\tau_{\ge n}$. Since $\tau_{\ge n+1}\subseteq\tau_{\ge n}$, there is a natural
transformation $P^n\to P^{n-1}$, which yields the \emph{slice tower} $\{P^nX\}_n$ of each
$G$-spectrum $X$.

Before proceeding, we introduce some terminology that will be useful in the computations that follow.
\begin{defn}\label{d:slice}\leavevmode
\begin{enumerate}
\item The homotopy fibre of $P^nX\to P^{n-1}X$, denoted $P^n_nX$, is the
\emph{$n$-slice} of $X$.
\item The homotopy fibre of $X\to P^{n-1}X$ is denoted $P_nX$.
\item We write $X\ge n$ if $X\in\tau_{\ge n}$, and $X\le n$ if $[C,X]^G=0$ for all
$C\in\tau_{\ge n+1}$ equivalently $[S^{k\rho_H+r},X]^H=0$ for all $H\le G$, $r\ge0$ and $k|H|>n$. 
\item We write $X\in[a,b]$ if $a\le X\le b$, and call $X$ an
\emph{$n$-slice} if $X\in[n,n]$.
\end{enumerate}
\end{defn}

Thus $P_nX\to X$ is the $\tau_{\ge n}$-colocalization and $X\to P^{n-1}X$ is the
$\tau_{\ge n}$-localization. Thus $P_nX\ge n$ and $P^{n-1}X\le n-1$. In particular
\[
X\ge n\iff P^{n-1}X\simeq*\iff P_nX\xrightarrow{\ \simeq\ }X,\quad
X\le n\iff X\xrightarrow{\ \simeq\ }P^nX .
\]

\begin{remark}\label{r:range}
If $X\in[a,b]$, then $P_aX\simeq X\simeq P^bX$, all slices outside $[a,b]$ vanish, and
the tower has finitely many stages:
\[
X=P^bX\to P^{b-1}X\to\cdots\to P^aX=P^a_aX\to P^{a-1}X=* .
\]
\end{remark}

We will repeatedly use the following results.
\begin{enumerate}[label=S\arabic*]
\item\label{S:HY} (Hill--Yarnall \cite[Theorem~A]{HY18}) A $G$-spectrum $X$ satisfies $X\ge n$ if and only if the non-equivariant homotopy groups $\pi_k\Phi^H(X)=0$ for all $H\le G$ and $k<n/|H|$. Here $\Phi^H(X)$ is the $H$-geometric fixed points of $i^*_HX$. Since $\Phi^H$ is symmetric monoidal with $\Phi^HS^V=S^{V^H}$, and $\Phi^H(H\underline M)$ is connective, we get: \emph{$S^V\wedge H\underline M\ge n$ whenever $\dim V^H\ge n/|H|$ for all $H$}, for every $V\in RO(G)$. In particular $\Sigma^mH\underline M\ge m$ for $m\ge0$.
\item\label{S:01} \cite{HHR16,U} $X$ is a $0$-slice if and only if $X\simeq H\underline M$; $\Sigma^1H\underline M$ is a $1$-slice if and only if all restrictions of $\underline M$ are injective; $\Sigma^{-1}H\underline M$ is a $(-n)$-slice if and only if transfers out of subgroups of order $\ge n$ are surjective and $\underline M(G/H)=0$ for $|H|<n$.
\item\label{S:rho} \cite[Corollary 4.25]{HHR16}If \(X\) is an \(n\)-slice, then \(\Sigma^\rho X\) is an \((n+|G|)\)-slice. i.e.
\[
P^{k+|G|}_{k+|G|}(\Sigma^{\rho_G}X)\simeq\Sigma^{\rho_G}P^k_k(X)
\]
where $\rho_G$ is the regular representation of group $G$.
\item  \cite[Proposition 4.26]{HHR16}]\label{lproduct}
$\tau_{\ge a}\wedge\tau_{\ge b}\subseteq\tau_{\ge a+b}$. In particular, if $X\ge0$ and $Y\ge m$, then $X\wedge Y\ge m$.
\item\label{S:res} \cite[Corollary 3.6]{Hil12} The restriction to a subgroup preserves slice towers; in particular, the restriction of a \(k\)-slice is again a \(k\)-slice.
\item \cite{U}\label{p:slices}
If $Y$ is a $k$-slice over $G/N$, then $\phi^*_NY$ is a $|N|k$-slice over $G$.
\item\label{S:D}\cite[Theorems.~I.7.7, I.7.8]{U}.
Let $I=I_{\mathbb Q/\mathbb Z}$ denote the Brown--Comenetz dual, i.e.\ the $G$-spectrum
representing $B\mapsto\mathrm{Hom}\big(\pi^G_0(X\wedge B),\mathbb Q/\mathbb Z\big)$, so that
$[B,IX]^G\cong\mathrm{Hom}\big(\pi^G_0(X\wedge B),\mathbb Q/\mathbb Z\big)$ and
$\underline\pi_j(IX)=\mathrm{Hom}(\underline\pi_{-j}X,\mathbb Q/\mathbb Z)$ with restriction
and transfer interchanged. Then, for every $G$-spectrum $X$ and every $n\in\mathbb Z$,
\[
X\le n\iff IX\ge -n,\qquad\qquad X\ge n\iff IX\le -n,
\]
and
\[
P^n_n(IX)\simeq I\,P^{-n}_{-n}X .
\]
In particular, the $(-n)$-slices are exactly the Brown--Comenetz duals of the $n$-slices.
\end{enumerate}
 
 In a cofibre sequence, a known slice bound for one spectrum can yield an equivalence between the other two spectra at a specified slice level.
\begin{lemma}\label{lPequiv}
Let $F\to X\to Y$ be a cofibre sequence of $G$-spectra and $m\in\mathbb Z$.
\begin{enumerate}
\item If $Y\le m-1$, then $P_mF\xrightarrow{\simeq}P_mX$. Hence
$P^k_kF\simeq P^k_kX$ for all $k\ge m$.
\item If $F\ge m+1$, then $P^mX\xrightarrow{\simeq}P^mY$. Hence
$P^k_kX\simeq P^k_kY$ for all $k\le m$.
\end{enumerate}
\end{lemma}
\begin{proof}
(1) Let $C\ge m$. Then $\Sigma C\ge m$ as well, so $[C,Y]=[\Sigma C,Y]=0$. The exact sequence
\[
\cdots[\Sigma C,Y]\to[C,F]\to[C,X]\to[C,Y]\cdots
\]
therefore gives $[C,F]\cong[C,X]$. Since
$P_m\to\mathrm{id}$ is the $\tau_{\ge m}$-colocalization, this yields natural bijections
$[C,P_mF]\cong[C,P_mX]$ for all $C\ge m$. As $P_mF$ and $P_mX$ both lie in $\tau_{\ge m}$, the
Yoneda lemma gives $P_mF\simeq P_mX$. For $k\ge m$ we have $P_k\simeq P_kP_m$, so
$P_kF\simeq P_kX$; taking cofibres of $P_{k+1}\to P_k$ gives the claim on slices.

\medskip
(2) Dually, let $Z\le m$. Then $[F,Z]=[\Sigma F,Z]=0$, so the exact sequence
\[
\cdots[\Sigma F,Z]\to[Y,Z]\to[X,Z]\to[F,Z]\cdots
\]
gives $[Y,Z]\cong[X,Z].$ Hence $P^mX\simeq P^mY$, and
$P^k\simeq P^kP^m$ for $k\le m$.
\end{proof}

\end{mysubsection}

\vspace{0.25cm}
Having established the necessary definitions and results for studying slices, we now determine bounds for a spectrum that identify the range in which its nontrivial slices can occur.
\medskip

\begin{mysubsection}{Slice ranges for $X(V,\uM)$}~
For a Mackey functor $\uM$ and a subgroup $H\leq G$, define
\[
b_H(\uM)
:=\min\bigl\{j\in\mathbb{Z}\mid
\pi_j\Phi^H H\uM\neq 0\bigr\}
\in\mathbb{Z}_{\geq 0}\cup\{\infty\},
\]
and call $b_H(\uM)$ the \emph{bottom degree} of the geometric fixed points. We set $b_H(\uM)=\infty$ if and only if $\Phi^H H\uM\simeq 0$. Let
\[
\uM^{\vee}:=\operatorname{Hom}(\uM,\mathbb{Q}/\mathbb{Z})
\]
denote the Pontryagin dual Mackey functor, obtained by interchanging restriction and transfer maps. We then define
\[
c_H(\uM):=b_H(\uM^{\vee}).
\]

\begin{remark}
Since
\[
\pi_0\Phi^H H\uM
\cong
\uM(H)\Big/
\sum_{K<H}\operatorname{tr}_K^H\uM(K),
\]
we have $b_H(\uM)=0$ if and only if
\[
\uM(H)\neq
\sum_{K<H}\operatorname{tr}_K^H\uM(K).
\]    
\end{remark}

We also introduce the following definitions.
\begin{myeq}{\label{boundn}}
n(V)=\min_{H\le G}\ |H|\big(d_H(V)+b_H(\uM)\big),
\end{myeq}
\begin{myeq}{\label{boundk}}
k(V)=\max_{H\le G}\ |H|\big(d_H(V)-c_H(\uM)\big),
\end{myeq}
the terms with $b_H=\infty$ (respectively \ $c_H=\infty$) being omitted. We use $d_H(V)$ to denote the dimension of the fixed point of $V$ for subgroup $H$, i.e. $ d_H(V)= \dim (V^H).$

Using the above notions, we establish bounds on the spectrum $X(V,\uM):=S^V\wedge H\uM$ and examine several particular cases.

\begin{thm}\label{trange}
Let $G$ be any finite group and $\uM$ be a Mackey functor. Then for any $V\in RO(G)$ with $X(V,\uM)=S^V\wedge H\uM\not\simeq *$, then we have
\begin{enumerate}
\item $X(V,\uM)\in [k(V),n(V)].$
\item Both the slices $P^{n(V)}_{n(V)}X(V,\uM)$ and $ P^{k(V)}_{k(V)}X(V,\uM)$ are non-zero.
\item If $\uM(e)\neq0$, this level is $|V|$, and $X(V,\uM)$ is a $|V|$-slice if and only if, for every $H\neq e$ with $b_H$ or $c_H$ finite,
\[
\Big\lceil\frac{|V|}{|H|}\Big\rceil-b_H(\uM)\ \le\ d_H(V)\ \le\ \Big\lfloor\frac{|V|}{|H|}\Big\rfloor+c_H(\uM).
\]
\end{enumerate}
\end{thm}
 
\begin{proof}
(1)  Since $\Phi^H$ is symmetric monoidal and $\Phi^H S^V\simeq S^{V^H}$, we have
\[
\Phi^H X(V,\uM)\simeq \Sigma^{d_H(V)}\Phi^H H\uM.
\]
This spectrum is connective, with its bottom nonzero homotopy group in degree
$d_H(V)+b_H(\uM)$. By~\eqref{S:HY}, the condition $X(V,\uM)\geq m$ holds if and only if
\[
d_H(V)+b_H(\uM)\geq
\left\lceil\frac{m}{|H|}\right\rceil
\]
for every subgroup $H\leq G$ with $b_H(\uM)<\infty$. Equivalently, this is the case if and only if
\[
m\leq |H|\bigl(d_H(V)+b_H(\uM)\bigr)
\]
for every subgroup $H$. 

 \medskip

 Since
\[
\underline{\pi}_j(IY)
\cong \operatorname{Hom}\bigl(\underline{\pi}_{-j}Y,\mathbb{Q}/\mathbb{Z}\bigr),
\]
with restriction and transfer maps interchanged, we have $IH\underline{M}\simeq H\underline{M}^{\vee}.$ 

Consequently,
\[
IX(V,\uM)
=F(S^V\wedge H\underline{M},I)
\simeq S^{-V}\wedge H\underline{M}^{\vee}.
\]
Applying above discussion to $\underline{M}^{\vee}$ and $-V$, we conclude that the bottom degree of $\Phi^H IX(V,\uM)$ is
\[
c_H(\underline{M})-d_H(V).
\]

By~\eqref{S:D}, the condition $X\leq m$ is equivalent to $IX\geq -m$.  This holds if and only if
\[
c_H(\underline{M})-d_H(V)
\geq \left\lceil-\frac{m}{|H|}\right\rceil
=-\left\lfloor\frac{m}{|H|}\right\rfloor
\]
for every subgroup $H\leq G$ with $c_H(\underline{M})<\infty$. Equivalently,
\[
m\geq |H|\bigl(d_H(V)-c_H(\underline{M})\bigr)
\]
for all such subgroups $H$, which is equivalent to $m\geq k(V)$, which proves (1).
\medskip

 (2) Put $n=n(V)$. Then $X\geq n$ but $X\not\geq n+1$. If $P_n^nX\simeq \ast$, the cofibre sequence
\[
P_{n+1}X\longrightarrow P_nX\longrightarrow P_n^nX
\]
implies that $P_nX\simeq P_{n+1}X$. Since $X\geq n$, we have $X\simeq P_nX$, and hence $X\simeq P_{n+1}X\geq n+1$, a contradiction.

\medskip

Similarly, suppose that $P_k^kX=0$, where $k=k(V)$. Then
\[
X\simeq P^kX\simeq P^{k-1}X\leq k-1,
\]
contradicting $X\not\leq k-1$. Hence both the slices are non-trivial, which completes the proof of~(2).

 \medskip
(3)  If $\uM(e)\neq0$, the term $H=e$ contributes $|V|$ to both (as $b_e=c_e=0$), so the common value is $|V|$. The displayed inequalities are $|V|\le|H|(d_H+b_H)$ and $|H|(d_H-c_H)\le|V|$ rewritten.
\end{proof}

Recall that a Mackey functor $\uM$ is called
\emph{cohomological} if
\[
\operatorname{tr}_K^H\circ\operatorname{res}_K^H
=|H:K|\cdot\operatorname{id}_{\underline{M}(H)}
\]
for all subgroups $K\leq H\leq G$. By Yoshida's theorem, cohomological
Mackey functors are precisely the modules over the constant Green functor
$\underline{\mathbb{Z}}$. Examples include
$\underline{\mathbb{Z}}$, $\underline{\mathbb{F}}_p$,
$\underline{\mathbb{Z}}^{*}$, and the Mackey functors
$\underline{g}_{\alpha}$.
 
\begin{cor}\label{ccohom}
Let $\uM$ be cohomological. Then $\Phi^HH\uM\simeq *$, so $b_H(\uM)=c_H(\uM)=\infty$, unless $|H|$ is a prime power. Consequently, in Theorem~\ref{trange} only subgroups of prime-power order contribute:
\[
n(V)=\min_{|H|=p^j}|H|\big(d_H(V)+b_H(\uM)\big),\qquad k(V)=\max_{|H|=p^j}|H|\big(d_H(V)-c_H(\uM)\big).
\]
\end{cor}

\begin{proof}
    By Yoshida's theorem, $H\uM$ is a module over $H\uZ$. Since $\Phi^H$ is strongly symmetric monoidal, $\Phi^HH\uM$ is a retract of $\Phi^HH\uZ\wedge\Phi^HH\uM$. Hence the result follows from~\cite[Proposition~11]{K21}.
\end{proof}

An immediate corollary for the $C_2$-case is as follows:
\begin{cor}\label{c2window}
  The $C_2$-spectrum  $S^{p+q\sigma}\wedge H\uF \in [\min(p+q,2p),\max(p+q,2p-4)]$. Moreover, it is a slice if and only if $0\le p-q\le 4.$
\end{cor}

We now apply Theorem~\ref{trange} to determine the slice range of
$X(V):=X(V,\uF)=\Sigma^V H\underline{\mathbb{F}}_2$ as a $K_4$-spectrum. Recall that
\[
\underline{\pi}_j(IY)
\cong
\operatorname{Hom}\bigl(\underline{\pi}_{-j}Y,\mathbb{Q}/\mathbb{Z}\bigr),
\]
with restriction and transfer maps interchanged. Consequently,
\[
IH\underline{\mathbb{F}}_2\simeq H\underline{\mathbb{F}}_2^{\,*},
\]
where $\underline{\mathbb{F}}_2^{\,*}$ is the Mackey functor that takes the
value $\mathbb{F}_2$ at every level, with all restriction maps equal to zero
and all transfer maps equal to the identity. By~\cite[Proposition~4.2]{GY20},
we have
\[
H\underline{\mathbb{F}}_2^{\,*}
\simeq
S^{3-\bar{\rho}}\wedge H\underline{\mathbb{F}}_2,
\qquad
\bar{\rho}=\alpha_0+\alpha_1+\alpha_2.
\]
It follows that
\[
\begin{aligned}
IX(V)
&=F\bigl(S^V\wedge H\underline{\mathbb{F}}_2,I\bigr)\\
&\simeq S^{-V}\wedge IH\underline{\mathbb{F}}_2\\
&\simeq S^{3-\bar{\rho}-V}\wedge H\underline{\mathbb{F}}_2.
\end{aligned}
\]

\begin{cor}\label{pk4range}
For $G=K_4$ and $V\in RO(K_4)$, all slices of
$X(V)=S^V\wedge H\underline{\mathbb{F}}_2$ lie in the range
$[n(V),k(V)]$, where
\begin{align*}
n(V)&=\min\bigl\{|V|,\ 2d_{H_i}(V),\ 4d_{K_4}(V)\colon i=0,1,2\bigr\},\\
k(V)&=\max\bigl\{|V|,\ 2d_{H_i}(V)-4,\ 4d_{K_4}(V)-12\colon i=0,1,2\bigr\}
\end{align*}
where $d_{H}=\dim (V^{H})$ for $H\subset K_4$. Moreover, $X(V)$ is a $|V|$-slice if and only if
$n(V)=k(V)=|V|$. Equivalently, this holds if and only if
\[
0\leq e_P\leq 4
\qquad\text{for all }P\in\{H_0,H_1,H_2\},
\]
where $e_P=2d_P-|V|$.
\end{cor}

  In particular, if $V$ be the trivial representation of dimension $n$, we obtain the following bounds.
\begin{cor}
   The $K_4$-spectrum $\Sigma^{n}H\uF$ is an $n$-slice if and only if $0\le n\le 4$. For $n\ge 4,$ we have $n\leq \Sigma^{n}H\uF \leq 4(n-3)$.
\end{cor}

\begin{proof}
    Note that, for $V=n,$ we have $e_H=n$ for all subgroups. Hence the result follows.
\end{proof}
Consider another example, for $V=n\alpha_i$ with $n\geq 0$. In this case, we obtain the following result.

 \begin{cor}{\label{bound0}}
   The $K_4$-spectrum $\Sigma^{n\alpha_i}H\uF$ is an $n$-slice if and only if $n=0$. For $n \ge 4,$ we have $0\leq \Sigma^{n\alpha_i}H\uF \leq 2n-4$ and $0\leq \Sigma^{n\alpha_i}H\uF \leq n$ if $0\le n \le 3.$
\end{cor}

\begin{proof}
We prove the result for $i=0$; the remaining cases follow by similar arguments. In this case, we have $e_{H_0}=n$, whereas $e_{H_1}=e_{H_2}=-n$. Therefore, using Corollary \ref{pk4range}, one derives $n=0$ is the only case when $X(n\alpha_0)$ is an $n$-slice.
\end{proof}
\end{mysubsection}

\section{Slices of Cofibres}

In~\cite[Proposition~5.9]{GK26}, the existence of the classes \(u_{\alpha_i}\) and \(v_i\), for \(i=0,1,2\), was established. For any \(V\in RO(K_4)\), we have four possible moves on \(X(V):=S^{V}\wedge H\uF\), obtained by multiplication or premultiplication by \(u_{\alpha_i}\) and \(v_i\)  (Proposition~\ref{pmoves}). In this section, we compute the fibres and cofibres of the maps induced by the multiplication of classes \(u_{\alpha_i}\) and \(v_{i}\), and subsequently determine the slices of the resulting spectra. 

\begin{mysubsection}{$u_{\alpha_i}$-twist}
~For $i=0,1,2$ and $V\in RO(K_4)$, define $X(V):=S^V\wedge H\underline{\mathbb{F}}_2,$ and let $C(V,\alpha_i)$ denote the cofibre of the map
\begin{myeq}{\label{mapu}}
    u_{\alpha_i}\colon X(V)\longrightarrow X(V+\alpha_i-1).
\end{myeq}

The following proposition proves that
cofibre $C(V\alpha_i)$ is homotopy equivalent to an Eilenberg--Mac Lane spectrum up to suspension. 

\begin{prop}\label{pCform}
The cofibre of the map in \eqref{mapu} is homotopy equivalent to the following spectrum
\[
C(V,\alpha_i)\simeq\Sigma^{V-1}H\ug_{(i)}.
\]
where $\ug_{(i)}$ is given in Example~\ref{Phimac}.

 Moreover, for $V=a+b_0\alpha_0+b_1\alpha_1+b_2\alpha_2$, we have
\begin{myeq}{\label{geopiece}}
C(V,\alpha_i)\simeq\mathrm{fib}\Big(\phi_{H_j}^*\big(\Sigma^{t+b_j\sigma}H_{C_2}\uF\big)\vee\phi_{H_k}^*\big(\Sigma^{t+b_k\sigma}H_{C_2}\uF\big)\longrightarrow\Sigma^tH\ug\Big)
\end{myeq}
where $t=a-1$, and $i,j,k$ are distinct elements of $[2]$.
\end{prop}
\begin{proof}
For convenience, we prove the result for $i=0$; the other cases follow similarly.

To understand the cofibre $C(V,\alpha_0)$, first observe that
\[
\upi_{k}(\Sigma^{1-\alpha_0}H\uF)\cong \begin{cases}
    \phi^{\ast}_{H_0}(\underline{f}_{C_2})\oplus \underline{f} & k=0\\
    0 & \text{otherwise}.
\end{cases}
\]
Then the homotopy groups of the cofibre $C(1-\alpha_0,\alpha_0)$ follow from the long exact sequence associated to the cofibre sequence 
\begin{myeq}\label{cofu0}
  \Sigma^{1-\alpha_0}H\uF \stackrel{u_{\alpha_0}}{\to} H\uF \to C(1-\alpha_0, \alpha_0)  
\end{myeq}
and, consequently, we obtain 
\[
\upi_{k}(C(1-\alpha_0,\alpha_0))\cong 
\begin{cases}
    \ug_{(0)} & k=0,\\
    0 & \text{otherwise}.
\end{cases}
\]
Smashing the following cofibre sequence \eqref{cofu0} with $S^{V+\alpha_0-1}$, we obtain
\[
C(V,\alpha_0)\simeq \Sigma^{V-1}H\ug_{(0)}.
\]
 For the second formula, consider the following short exact sequence of Mackey functors:
\[
0\longrightarrow \ug_{(0)}
\longrightarrow \phi_{H_1}^*\uF
\oplus \phi_{H_2}^*\uF
\longrightarrow \ug
\longrightarrow 0.
\]
Smashing the associated cofibre sequence of Eilenberg--Mac Lane spectra with
$S^W$, where
$
W=t+b_1\alpha_1+b_2\alpha_{2}$ and using the equivalences
\[
S^W\wedge H\phi_{H_1}^*\uF
\simeq
\phi_{H_1}^*\bigl(S^{W^{H_1}}\wedge H\uF\bigr),
\]
where $W^{H_1}=t+b_1\sigma$ with $\sigma$ denoting the sign representation of $K_4/{H_1}$, we obtain the desired expression. Also, using the analogous equivalence for ${H_2}$ along with the following equivalence $S^W\wedge H\ug
\simeq \Sigma^t H\ug.$
\end{proof}

In the following proposition, we characterize when the spectrum $\Sigma^{V}H\ug_{(i)}$ satisfy a given \(m\)-slice bound. We then establish equivalences between the corresponding slices as a consequence of Lemma~\ref{lPequiv}.

\begin{prop}\label{psusgi}
Let $V \in RO(K_4)$. 
\begin{enumerate}
    \item The spectrum \(\Sigma^V H\ug_{(i)}\leq m\)  if and only if
\[
\max (2d_{H_j},\;2d_{H_k},\;4(d_{K_4}-1))\leq m.
\]
Moreover, if it holds, then the map \(\ref{raiseu}\), $
u_{\alpha_i}\colon X(V+1-\alpha_i)\longrightarrow X(V)$ induces equivalences $P_{m+1}X(V+1-\alpha_i)\simeq P_{m+1}X(V)$
and $
P^{m+1}_{m+1}X(V+1-\alpha_i)
\simeq P^{m+1}_{m+1}X(V).$

\item The spectrum \(\Sigma^{V-2}H\ug_{(i)}\) is greater than or equal to \(m\) if and only if 
\[\min(2(d_{H_j}-2),
2(d_{H_k}-2),
4(d_{K_4}-2))\ge m.\]
Moreover, if it holds, then the map \(\ref{loweru}\), $u_{\alpha_i}:X(V)\xrightarrow{\ \ }X(V+\alpha_i-1)$ induces equivalences $
P^{m+1}X(V)\simeq P^{m+1}X(V+\alpha_i-1)$ and $
P^{m+1}_{m+1}X(V)
\simeq P^{m+1}_{m+1}X(V+\alpha_i-1).$
\end{enumerate} 
\end{prop}

\begin{proof}
(1) Consider the cofibre sequence
\[
\Sigma^{1-\alpha_i}H\uF\xrightarrow{u_{\alpha_i}}H\uF\to H\ug_{(i)}.
\]
The identification $IH\uF\simeq\Sigma^{2-\alpha_j-\alpha_k}\Sigma^{1-\alpha_i}H\uF$, gives the fibre sequence
\[
IH\ug_{(i)}\to\Sigma^{2-\alpha_j-\alpha_k}\Sigma^{1-\alpha_i}H\uF\xrightarrow{\Sigma^{2-\alpha_j-\alpha_k}u_{\alpha_i}}\Sigma^{2-\alpha_j-\alpha_k}H\uF ,
\]
since $Iu_{\alpha_i}$ is the nonzero element of $\pi^{K_4}_{1-\alpha_i}H\uF\cong\bF\{u_{\alpha_i}\}$. Hence
\[
I(\Sigma^VH\ug_{(i)})\simeq\Sigma^{1-\alpha_j-\alpha_k-V}H\ug_{(i)}.
\]

 If $\Sigma^V H\ug_{(i)} \le m$, then $\Sigma^{1-\alpha_j-\alpha_k-V}H\ug_{(i)} \ge -m.$ Thus $\pi_q\Phi^J\Sigma^{1-\alpha_j-\alpha_k-V}H\ug_{(i)} =0$ for $q<\lceil(-m)/|J|\rceil$. Therefore the conditions
read $-d_{H_j}\ge-\lfloor m/2\rfloor$, $-d_{H_k}\ge-\lfloor m/2\rfloor$ and $1-d_{K_4}\ge-\lfloor m/4\rfloor$. Hence the result follows.

(2) Note that  $\Phi^J \Sigma^{V-2} H\ug_{(i)}\simeq\Sigma^{d_J-2}\Phi^JH\ug_{(i)}.$ If $J\in\{e,H_i\}$, then $\res^{K_4}_J\ug_{(i)}=0$, so $\Phi^J \Sigma^{V-2} H\ug_{(i)}=0$ and there is no condition. If $J\in\{H_j,H_k,K_4\}$, then
$\Phi^JH\ug_{(i)}$ is connective. So the bottom homotopy group of $\Phi^J\Sigma^{V-2} H\ug_{(i)}$ is in degree exactly $d_J-2$, and
the condition at $J$ is $d_J-2\ge\lceil m/|J|\rceil$, i.e.\ $|J|(d_J-2)\ge m$. For $J=H_j,H_k,K_4$ this gives the three inequalities. 
\end{proof}

Since the restriction to \({H_i}\) of the cofibre of the map~\eqref{mapu}, namely
\[
\Sigma^{V-1}H\ug_{(i)}
\simeq \Sigma^{t+\Sigma_{i\neq \ell=0}^{2}b_\ell\alpha_\ell}H\ug_{(i)},
\]
is trivial, while its restrictions to \({H_j}\) and \({H_k}\) are \(\Sigma^{t+b_j}H_{C_2}\ug_{C_2}\) and \(\Sigma^{t+b_k}H_{C_2}\ug_{C_2}\), respectively, these yield \(2(t+b_j)\)- and \(2(t+b_k)\)-slices, respectively. All the remaining slices are geometric; hence, the dimensions of the remaining slices are given by multiples of \(|K_4|=4\).

To establish the main theorem of this section, we need to compute the \(C_2\)-slices of the spectra \(\Sigma^{a+b\sigma}H_{C_2}\uF\). However, by \(\rho_{C_2}\)-periodicity, it suffices to compute the slices of spectra of the form \(\Sigma^{a-b}H_{C_2}\uF\).
\begin{prop}{\label{C_2slices}}
For \(n\in\mathbb Z\), the following hold:
\begin{enumerate}
    \item The slices of \(\Sigma^n H_{C_2}\uF\) are concentrated in the range $[\min(n,2n),\,\max(n,2n-4)]$. 
    \item For \(n\geq 4\), the nontrivial slices are given by
    \[
    P_k^k(\Sigma^n H_{C_2}\uF)\simeq
    \begin{cases}
    \Sigma^{(\frac n2-2)\rho_{C_2}+4}H_{C_2}\uF
        & k=n,\quad n\text{ even},\\
    \Sigma^{(\frac{n-3}{2})\rho_{C_2}+3}H_{C_2}\uF
        & k=n,\quad n\text{ odd},\\
    \Sigma^{\frac k2}H_{C_2}\ug_{C_2}
        & k=n+i,\quad
        \begin{array}{l}
        i=2,4,\ldots,n-4,\quad n\text{ even},\\
        i=1,3,\ldots,n-4,\quad n\text{ odd},
        \end{array}\\
    *
        & \text{otherwise}.
    \end{cases}
    \]

    \item For \(n<0\), the nontrivial slices are given by
    \[
    P_k^k(\Sigma^n H_{C_2}\uF)\simeq
    \begin{cases}
    \Sigma^{\frac n2\rho_{C_2}}H_{C_2}\uF
        & k=n,\quad n\text{ even},\\
    \Sigma^{1+\frac{n-1}{2}\rho_{C_2}}H_{C_2}\uF
        & k=n,\quad n\text{ odd},\\
    \Sigma^{\frac k2}H_{C_2}\ug_{C_2}
        & k=n-i,\quad
        \begin{array}{l}
        i=2,4,\ldots,-n,\quad n\text{ even},\\
        i=1,3,\ldots,-n,\quad n\text{ odd},
        \end{array}\\
    *
        & \text{otherwise}.
    \end{cases}
    \]
\end{enumerate}
Moreover, the spectra $\Sigma^{a+b\sigma}H_{C_2}\uF$ is $(a+b)$-slice if and only if $b\leq a\leq b+4$.
\end{prop}
\begin{proof}
\begin{enumerate}
    \item  The bounds for the \(C_2\)-spectra follow directly from Corollary~\ref{c2window}.

\item The result for \(n\geq 4\) follows from~\cite[Theorem~3.18]{GY20}.

\item For \(n>0\), it follows from~\eqref{S:D} that
\[
\begin{aligned}
P_{-n}^{-n}(\Sigma^{-n}H_{C_2}\uF)
&\simeq P_{-n}^{-n}\bigl(I(\Sigma^{n+2-2\sigma}H_{C_2}\uF)\bigr)\\
&\simeq I\bigl(\Sigma^{-2\rho_{C_2}}
P_{n+4}^{n+4}(\Sigma^{n+4}H_{C_2}\uF)\bigr).
\end{aligned}
\]
When \(n\) is even, we have $P_{n+4}^{n+4}(\Sigma^{n+4}H_{C_2}\uF)
\simeq \Sigma^{\frac n2\rho_{C_2}+4}H_{C_2}\uF$, which gives $P_{-n}^{-n}(\Sigma^{-n}H_{C_2}\uF)
\simeq \Sigma^{-\frac n2\rho_{C_2}}H_{C_2}\uF.$

On the other hand, when \(n\) is odd, we have $P_{n+4}^{n+4}(\Sigma^{n+4}H_{C_2}\uF)
\simeq \Sigma^{\frac{n+1}{2}\rho_{C_2}+3}H_{C_2}\uF$, which gives $P_{-n}^{-n}(\Sigma^{-n}H_{C_2}\uF)
\simeq \Sigma^{1-\frac{n+1}{2}\rho_{C_2}}H_{C_2}\uF.$

\end{enumerate}
This completes the proof for $n$-slice. The other remaining cases follow similarly. 
\end{proof}
In the following theorem, we determine the slices of the cofibre $C(V,\alpha_i)$. 
\begin{thm}\label{tCslices} Let $V=a+b_0\aO+b_1\aI+b_2\aII$ and $t:=a-1$. 
\begin{enumerate}
    \item The slices of $C(V,\alpha_i)$ for every $n\notin [4t-4,4t]$,  
\[
P^n_nC(V,\alpha_i)\simeq
\begin{cases}
\phi^*_{H_j}\big(P^{n/2}_{n/2}Y_{H_j}\big)\vee\phi^*_{H_k}\big(P^{n/2}_{n/2}Y_{H_k}\big) & \text{n is even} \\
\ast & \text{n is odd}
\end{cases}
\]
where $Y_{H_j}=\Sigma^{t+b_j\sigma}H_{C_2}\uF$ and $Y_{H_k}=\Sigma^{t+b_k\sigma}H_{C_2}\uF$. 
\item The slices of \(C(V,\alpha_i)\) in degrees \(4t-4\), \(4t-2\), and \(4t\) are described as follows. Assume, without loss of generality, that \(b_j\geq b_k\), and set \(e=t-b_j\).
\begin{enumerate}
    \item For $b_j>t$,
\[
P^{4t}_{4t}C(V,\alpha_i)\simeq
\begin{cases}\Sigma^tH\ug & b_k>t\\ \Sigma^{t\rho}H\phi^*_{H_k}\uF & b_k=t\\ * & b_k<t\end{cases}
\]
and for 
$\ell\in\{4t-4,4t-2\}$,

\[
P^\ell_\ell C(V,\alpha_i)\simeq \begin{cases}
\Sigma^{t+b_k\alpha_k} H\phi^*_{H_k}\uF & \ell=2(t+b_k)\\
\ast & \text{otherwise.}
\end{cases}
\]

\item Assume that  $b_j=t-e$ with $e\in\{0,1,2\}.$ In level $4t-2e$,
\[
P^{4t-2e}_{4t-2e}C(V,\alpha_i)\simeq
\begin{cases}S^{(t-e)\rho}\wedge\Sigma^eH\phi^{\ast}_{H_j}\uf_{C_2} & b_k<b_j\\ S^{(t-e)\rho}\wedge\Sigma^eH\ug_{(i)} & b_k=b_j.\end{cases}
\]
If $t-2\le b_k<b_j$, the level $2(t+b_k)$ carries $S^{t+b_k\alpha_k}\wedge H\phi^*_{H_k}\uF$. All other band levels are $*$.

\item For $b_j\le t-3$,  $P^{4t-4}_{4t-4}C(V,\alpha_i)\simeq\Sigma^{t-1}H\ug$, and $P^{4t-2}_{4t-2}C(V,\alpha_i)\simeq P^{4t}_{4t}C(V,\alpha_i)\simeq*$.
\end{enumerate}
\end{enumerate}
\end{thm}
\begin{proof}Again, we prove it for $i=0$; other cases follow in a similar manner. 

(1) Using Proposition \ref{pCform}, we have a cofibre sequence \[
C(V,\aO)\to \phi^*_{H_1}Y_{H_1}\vee\phi^*_{H_2}Y_{H_2} \to\Sigma^tH\ug.
\]
Since the cofibre $\Sigma^tH\ug$ is a $4t$-slice, so $P_nC(V,\aO)\simeq P_nB$ for $n>4t$.
Furthermore, the fibre
$\Sigma^{t-1}H\ug$ is a $4(t-1)$-slice, so $P^nC(V,\aO)\simeq P^nB$ for $n<4t-4$. Therefore, we have
\[
P^n_nC(V,\aO)\simeq\phi^*_{H_1}\big(P^{n/2}_{n/2}Y_{H_1}\big)\vee\phi^*_{H_2}\big(P^{n/2}_{n/2}Y_{H_2}\big)\qquad\text{for } n\notin[4t-4,4t],
\]
and $P^n_nC(V,\aO)=0$ for $n$ odd.

(2) 
Put $c_\nu=b_\nu-t+1$. The class $a_{\aO}$ is invertible on $H\ug_{(0)}$, so
\[
Z:=S^{-(t-1)\rho}\wedge C(V,\aO)\simeq S^{1+c_1\aI+c_2\aII}\wedge H\ug_{(0)} .
\]
By $\rho$-periodicity~\eqref{S:rho}, the slices of $C(V,\aO)$ in levels $4t-4,4t-2,4t$ are $S^{(t-1)\rho}\wedge P^m_mZ$ for $m=0,2,4$. The three cases of the
theorem become $c_j\ge2$, $c_j=1-e$ with $e\in\{0,1,2\}$, and $c_j\le-2$.
Smash $0\to\phi_{H_{j}}^{\ast}\uf_{C_2}\to\ug_{(0)}\to\phi^*_{H_k}\uF\to0$ with $S^{1+c_1\aI+c_2\aII}$, we derive a cofibre sequence

\[
\phi^*_{H_j}\big(S^{2+(c_j-1)\sigma}\wedge H_{C_2}\uF\big)= :F\to Z \to  E:= \phi^*_{H_k}\big(S^{1+c_k\sigma}\wedge H_{C_2}\uF\big).
\]

{\noindent \textbf{Case I:}  For $c_j\ge 2$.

By Corollary \ref{c2window}, the $C_2$-spectrum $S^{2+(c_j-1)\sigma}\wedge H_{C_2}\uF \ge \min(c_j+1, 4)=3$ which give $F \ge 6.$ Hence, by Lemma~\ref{lPequiv}
\[
P^m_m Z\simeq P^m_m E \text{ for } m \le 4.
\]
Now we will analyse the spectrum $E$.
\begin{enumerate}[label=\textup{(\roman*)}]
    \item If $c_k\ge 2,$ then by Corollary~\ref{c2window}, $E\ge 2\min(1+c_k, 2)\ge 4.$ The $2$-slice of $S^{1+c_k\sigma}\wedge H_{C_2}\uF$ is given by $H_{C_2}\ug_{C_2}.$ Hence the $4$-slice of $E$ is $ H\ug$.

    \item If $c_k=1,$ then the $2$-slice of $S^{1+\sigma}\wedge H_{C_2}\uF$ is $S^{\rho_{C_2}}\wedge H_{C_2}\uF$, then the $4$-slice of $Z$ is $S^{\rho}\wedge H\phi^*_{H_k}\uF.$

    \item If $c_k\in \{-1, 0\}$, then $E$ is a single slice in dimension $2(1+c_k).$ One computes this slice is $S^{1+c_k\alpha_k}\wedge H\phi^*_{H_k}\uF$.

    \item If $c_k\le -2$, then $E\le2\max(1+c_k,-2)\le-2$. Hence $E$ does not have slices in the dimensions $0, 2,4.$ Therefore, $P^m_m Z\simeq \ast$ for $m=0,2,4.$
\end{enumerate}

\noindent \textbf{Case II:}  $c_j=1-e$, with $e\in \{0,1,2\}.$

In this case, $F\simeq \phi^*_{H_j}(S^{2-e\sigma}\wedge H_{C_2}\uF).$ By Corollary \ref{c2window}, it is a $(4-2e)$-slice. Thus we note that $F\simeq S^{(1-e)\rho}\wedge\Sigma^eH\phi^{\ast}_{H_j}\uf_{C_2}$. To understand the slices of $Z$, we do case-by-case analyzes of $E$.

\begin{enumerate}[label=\textup{(\roman*)}]
    \item If $c_k<c_j$, then $E \le2\max(1+c_k,-2)\le2-2e<4-2e$. Therefore,  $F=P_{\ge4-2e}Z$ and $E=P^{3-2e}Z$. The band slices of $Z$ are $F$
in level $4-2e$, together with $E=S^{1+c_k\alpha_k}\wedge H\phi^*_{H_k}\uF$ in level $2(1+c_k)$ when $c_k\in\{-1,0\}$. For $c_k\le-2$, $E$
lies below level $0$.
\item If $c_k=c_j$, then $E=\phi^*_{H_k}(S^{1+(1-e)\sigma}\wedge H_{C_2}\uF)$ is a $(4-2e)$-slice. Hence $Z$ is also a $(4-2e)$-slice, and
\[
Z\simeq S^{1+(1-e)(\aI+\aII)}\wedge H\ug_{(0)}\simeq S^{(1-e)\rho}\wedge\Sigma^eH\ug_{(0)} .\]
\end{enumerate}

\noindent \textbf{Case III:} $c_j\leq -2$.

For this case, we have $E\le 2\max(1+c_k, -2)=-2.$ Hence $P^m_mZ \simeq P^m_m F$ for $m \ge 0.$ Note that $F \le 2\max(c_j+1, 0)=0$. Therefore $P^2_2 Z\simeq P^4_4 Z\simeq \ast.$ Thus we need to determine the $0$-slice of $F$. This will be equivalent to computing the $0$-slice of the $C_2$-spectrum $S^{2+(c_j-1)\sigma}\wedge H_{C_2}\uF,$ which is given by  $H_{C_2}\upi_0(S^{2+(c_j-1)\sigma}\wedge H_{C_2}\uF)$. Note that $\upi_0S^{2+(c_j-1)\sigma}\wedge H_{C_2}\uF$ is $0$ at $C_2/e$ (since the underlying dimension $1+c_j\ne 0$) and
$\F\langle\Sigma^{-1}\frac1{a^{-c_j}u}\rangle$ at $C_2/C_2$, that is, $\ug_{C_2}$. Hence $P^0_0 Z\simeq H\phi^*_{H_k}\ug_{C_2}\simeq H\ug.$

Now we smash $S^{(t-1)\rho}$ with $Z$ to get back the original slices of $C(V, \alpha_0)$. 

Hence the result follows.}
\end{proof}
\end{mysubsection}

\begin{example}
    Consider $X(5+2\alpha_1+2\alpha_2)=S^{5+2\alpha_1+2\alpha_2}\wedge H\uF$, then its slice tower is
    \[
\xymatrix{
P^{10}_{10}\simeq \Sigma^{3+2\alpha_1+2\alpha_2} H\ug_{(0)} \ar[r]& P^{10}X(5+2\alpha_1+2\alpha_2)\simeq X(5+2\alpha_1+2\alpha_2)\ar[d]^{u_{\alpha_0}} \\ & P^{9}X(5+2\alpha_1+2\alpha_2)\simeq X(4+\alpha_0+2\alpha_1+2\alpha_2).
}
\]
\end{example}

\begin{mysubsection}{$v_i$-twists} ~For $i\in[2]$, let $\tau_i=1+\alpha_i-\alpha_j-\alpha_k\in RO(K_4),$ where $\{j,k\}=[2]\setminus\{i\}$. Consider the following map of $H\underline{\mathbb{F}}_2$-modules:
\begin{myeq}\label{dcvi}
v_i\colon \Sigma^{\tau_i}H\underline{\mathbb{F}}_2
\longrightarrow H\underline{\mathbb{F}}_2,
\end{myeq}
and let $C(v_i)$ denote its cofibre. One can compute from~\cite[Theorem~5.19]{GK26}, at $K_4/K_4$-level, $\pi_{0}(\Sigma^{\tau_{i}}H\uF)\cong \F\{v_{j}\Sigma^{-1}\frac{1}{a_{\alpha_k}u_{\alpha_k}}\}$, whose restriction to the different subgroup is defined as follows.

\[
\res_{H_\ell}(v_{j}\Sigma^{-1}\frac{1}{a_{\alpha_k}u_{\alpha_k}}) =\begin{cases}
    \Sigma^{-1}\frac{1}{a_{\sigma}u_{\sigma}} & \ell=i,\\
    0 & \ell=j,k.
    \end{cases}
\]
Hence, $\Sigma^{\tau_i}H\uF\simeq H\uN_{H_i}$ where $\underline{N}_{H_i}$ is described in Example~\ref{MacN}. 

On $\underline\pi_0$ the map $v_i$ is $\varphi\colon\uN_{H_i}\to\uF$, with
$\ker\varphi\cong\phi^*_{H_i} {\uF}$ and $\coker \varphi \cong \phi^{*}_{H_i}\uF$, because the map $\phi$ acts nontrivially on orbits $K_4/e$, $K_4/H_j$ and $K_4/H_k$, whereas acts trivially on $K_4/H_i$ and $K_4/K_4$ orbits.

\begin{prop}\label{pcvi}
The cofibre $C(v_i)$ has exactly two nonzero slices,
\[
P^0_0\,C(v_i)\simeq H\phi^{*}_{H_i}\uF ,\qquad
P^2_2\,C(v_i)\simeq \Sigma H\phi^*_{H_i}\uF\simeq\phi^*_{{H_i}}\big(\Sigma H_{C_2}\uF\big),
\]
and $P^\ell_\ell\,C(v_i)\simeq *$ for $\ell\neq0,2$. The slice tower is the cofibre sequence
\begin{myeq}\label{extcvi}
    \Sigma H\phi^*_{H_i}\uF\longrightarrow C(v_i)\longrightarrow H\phi^*_{H_i}\uF .
   \end{myeq} 
\end{prop}

\begin{proof} Again, we do the proof for $i=0$, other cases follow in similar manner.

Both $\Sigma^{\tau_0}H\uF\simeq H\underline{N}_{H_0}$ and $H\uF$ are Eilenberg--Mac Lane, so the long exact sequence associated from the map~\eqref{dcvi} for $i=0$ gives
$\underline\pi_0C(v_0)=\coker\varphi$, $\underline\pi_1C(v_0)=\ker\varphi$, and
$\underline\pi_mC(v_0)=0$ for $m\neq0,1$. 

By \eqref{p:slices}, Geometric inflation along $K_4\to K_4/H_0\cong C_2$ multiplies slice levels by $|H_0|=2$ . Over $C_2$, $H_{C_2}\uF$ is 
$0$-slice. Using~\eqref{S:01}, $\Sigma H_{C_2}\uF$ is a $1$-slice. Hence, $H\phi^*_{H_0}\uF$ and $\Sigma H\phi^*_{H_0}\uF$ are $0$-slice and $2$-slice respectively, as $K_4$ spectrum.

In a cofibre sequence $A\to C\to B$ with $A\ge 2$ and $B\le 1$, we have $P_2C\simeq A$ and $P^1C\simeq B$.
Since $B$ is a $0$-slice and $A$ a $2$-slice, the claim follows.
\end{proof}

 Since the cofibre of the map
\[
\Sigma^{V+\tau_i}H\underline{\mathbb{F}}_2
\longrightarrow
\Sigma^V H\uF
\]
is $S^V\wedge C(v_i)$, the following corollary describes this spectrum and determines its slices. 

\begin{cor}\label{csphcvi} For $V\in RO(K_4)$, then $S^V\wedge C(v_i)$ is an extension
\[
\phi^*_{H_i}\big(S^{\,1+V^{H_i}}\wedge H_{C_2}\uF\big)\longrightarrow S^V\wedge C(v_i)\longrightarrow
\phi^*_{H_i}\big(S^{\,V^{H_i}}\wedge H_{C_2}\uF\big)
\]
where $V^{H_i}$ is regarded as a $K_4/H_i\cong C_2$-representation.
\end{cor}

 In the following proposition, we obtain conditions for the map induced by multiplication by \(v_i\).  
\begin{prop}\label{psuscvi}
Let $V \in RO(K_4)$, then 

\begin{enumerate}
    \item  $\Sigma^V C(v_i) \le m$ if and only if
    \[
 2(1+d_{H_i})\le m\text{ and } 4(d_{K_4}-1)\le m . 
    \]
Moreover, if it holds, then the map \ref{raisev}, $v_{\alpha_i}:X(V+\tau_i) \to X(V)$ induces isomorphisms $P_{m+1}X(V+\tau_i)\simeq P_{m+1}X(V)$ and  $P^{m+1}_{m+1}X(V+\tau_i)\simeq P^{m+1}_{m+1}X(V).$

\item $\Sigma^{V-\tau_i-1} C(v_i) \ge m$ if and only if 
\[
2d_{H_i}\ge  m+6\text{ and }4d_{K_4} \ge m+8.
\]
Moreover, if it holds, then the map \ref{lowerv}, $v_i:X(V)\xrightarrow{\ \ }X(V-\tau_i)$ induces isomorphism $P^{m+1} X(V) \simeq P^{m+1}X(V-\tau_i)$ and $P^{m+1}_{m+1}X(V)\simeq P^{m+1}_{m+1}X(V-\tau_i).$

\item For any $V\in RO(K_4)$, $\Sigma^{V}C(v_i)$ can never be a slice.
\end{enumerate}
\end{prop}
\begin{remark}
    Note that \(\res_{H_j}C(v_i)\simeq \res_{H_k}C(v_i)\simeq *\). Hence, it suffices to verify the following proposition at the \(H_i\)- and \(K_4\)-levels.
\end{remark}
\begin{proof}[Proof of Proposition~\ref{psuscvi}]
   (1) Using \eqref{dcvi}, we have a fibre sequence
\[
IC(v_i)\to\Sigma^{3-\bar\rho}H\uF\xrightarrow{\ \Sigma^{3-\bar\rho-\tau_i}v_i\ }\Sigma^{3-\bar\rho-\tau_i}H\uF ,
\]
since $Iv_i$ is the nonzero element of $\pi^{K_4}_{\tau_i}H\uF\cong\bF$. Hence
\[
IC(v_i)\simeq\Sigma^{2-\bar\rho-\tau_i}C(v_i)=\Sigma^{1-2\alpha_i}C(v_i)
\mbox{ and } I\Sigma^V C(v_i)\simeq \Sigma^{1-2\alpha_i-V}C(v_i).\]

\noindent By \cite{U}, $\Sigma^V C(v_i)\le m$ if and only if $\Sigma^{1-2\alpha_i-V}C(v_i) \ge -m.$ Thus $\pi_q\Phi^J\Sigma^{1-2\alpha_i-V}C(v_i)=0$ for $q< \lceil -m/|J| \rceil.$ The representation $1-2\alpha_i-V$ has fixed dimensions $-1-d_{H_i}$ at $H_i$
and $1-d_{K_4}$ at $K_4$. By the lower bound, the condition is $2(1+d_{H_i})\le m$ and $4(d_{K_4}-1)\le m$.

(2) If $\Sigma^{V-\tau_i-1}C(v_i)\geq m$, then we have $\pi_{q}\Phi^{J}\Sigma^{V-\tau_i-1}C(v_i)=0$ for $q<\lceil m/|J|\rceil$. The representation $V-\tau_i-1$ has fixed point $d_{H_i}-2$ at $H_i$ and $d_{K_4-2}$ at $K_4$. By the lower bound, the condition $2d_{H_i}\geq m+4$ and $4d_{K_4}\geq m+8$. 

(3) For $\Sigma^V C(v_i)$ to be an $m$-slice, both its upper and lower slice bounds must equal $m$. However, by (2), the condition $\Sigma^{(V+\tau_i+1)-\tau_i-1} C(v_i)\geq m$ is equivalent to $m\leq \min(2d_{H_i},4d_{K_4}).$ On the other hand, (1) implies that
\[
m\geq 2(1+d_{H_i}).
\]
These inequalities are incompatible. Hence, $\Sigma^V C(v_i)$ cannot be an $m$-slice.
\end{proof}

\begin{example}
Consider the spectrum
\[
X(4+\alpha_0)=S^{4+\alpha_0}\wedge H\uF.
\]
By Corollary~\ref{pk4range}, this spectrum lies in the slice interval $[5,6]$. We consider the following diagram of spectra, in which each pair of consecutive horizontal and vertical maps forms a cofibre sequence:
\[
\xymatrix{
 \Sigma^{3} H\phi_{H_0}^{\ast}\uF\ar[r]&\Sigma^{2+\alpha_1+\alpha_2}C(v_{0})\ar[d]  \ar[r]& P^{6}X(4+\alpha_0)\simeq X(4+\alpha_0)\ar[d]^{v_0} \\
 & \Sigma^{2}H\phi_{H_0}^{\ast}\uF & X(3+\alpha_1+\alpha_2).
}
\]
Here, $\Sigma^{3} H\phi_{H_0}^{\ast}\uF$, $\Sigma^{2} H\phi_{H_0}^{\ast}\uF$ are $6$-slice and $4$-slice respectively, whereas $X(3+\alpha_1+\alpha_2)$ is $5$-slice. However, this is not a slice tower for $X(4+\alpha_0)$.

To determine $5$ and $6$-slice, note that $ P_{k}^{k}(X(4+\alpha_0))\simeq \Sigma^{\rho}P^{k-4}_{k-4}(X(3-\alpha_1-\alpha_2))$ for $k=5,6$. So it is sufficient to compute $1$- and $2$-slice of $X(3-\alpha_1-\alpha_2)\in [1,2]$. One can compute these slices by restricting them to order $2$-subgroups, using the cofibre sequence $P^n_n\to P^n \to P^{n-1}$ and properties of $1$-slice~\eqref{S:01}.

Hence, the slice tower for $X(4+\alpha_0)$ is
\[
P^{6}_{6}\simeq\xymatrix{
\Sigma^{1+\rho}H\phi_{H_0}^{\ast}\uf_{C_2}  \ar[r]& P^{6}X(4+\alpha_0)\simeq X(4+\alpha_0)\ar[d] \\
 & P^{5}_{5}\simeq \Sigma^{1+\rho}H\uf.
}
\]
\end{example}
\end{mysubsection}

\section{Slices of $X(\infty \alpha_i)$}

Consider \(F(\alpha_i;n)\) to be the fibre of the map \(X(n\alpha_i)\to X(\infty\alpha_i)\), induced by the inclusion \(S^{n\alpha_i}\hookrightarrow S^{\infty\alpha_i}\), or equivalently, by iterated multiplication by \(a_{\alpha_i}\). Here, \(X(n\alpha_i)=S^{n\alpha_i}\wedge H\uF\). Thus, we have the fibre sequence
\begin{myeq}\label{fibalpha_0}
F(\alpha_i;n)\longrightarrow X(n\alpha_i)\longrightarrow X(\infty\alpha_i).
\end{myeq}

We prove each statement for \(i=0\); the cases \(i=1,2\) follow by applying automorphisms of \(K_4\), which permute the subgroups \((H_0,H_1,H_2)\) and the corresponding representations \((\alpha_0,\alpha_1,\alpha_2)\).

\begin{prop}\label{prop:Xinf}
For $i=0,1,2$, let $\mathcal{F}_i = \{e, H_i\}$. Then:
  \begin{enumerate}
    \item There are equivalences of $K_4$-spectra
      \[
        X(\infty\alpha_i) \simeq \EFt_{i} \wedge H\uF
        \simeq \Phi_{H_i}^*\bigl(\widetilde{EC_2} \wedge H_{C_2}\uF\bigr),
      \]
      and its homotopy Mackey functors are
      \[
        \upi_k X(\infty\alpha_i) \cong
        \begin{cases}
          \ug_{(i)} & k \geq 0,\\
          0    & k < 0;
        \end{cases}
      \]
    \item For $n\geq 0$, there is an equivalence
      \[
        F(\alpha_i; n) \simeq E{\mathcal{F}_{i}}_+ \wedge S^{n\alpha_i} \wedge H\uF
        \simeq \Phi_{H_i}^*\bigl({EC_2}_+ \wedge S^{n\sigma} \wedge H_{C_2}\uF\bigr).
      \]
  \item  For $n\geq 0$, the homotopy Mackey functors of the fibre $F(\alpha_i; n)$ are
  \[
    \upi_k F(\alpha_i; n) \cong
    \begin{cases}
      0                     & k < n,\\
      \ug_{(i)}\oplus\Phi_{H_i}^*\underline{f}_{C_2} & k = n,\\
      \ug_{(i)}                  & k > n.
    \end{cases}
  \]\end{enumerate}
\end{prop}
\begin{proof}Proof for $i=0$.

\begin{enumerate}
    \item Since $X(\infty\aO) \simeq S^{\infty\aO} \wedge H\uF$, the equivalence
follows from~\eqref{PhiL}.

For homotopy Mackey functors, we obtain the following expression from~\cite{BG21}.
\[
  \upi_k\bigl(\widetilde{EC_2} \wedge H_{C_2}\uF\bigr) \cong
  \begin{cases}
    \ug_{C_2} & k \geq 0,\\
    0  & k < 0.
  \end{cases}
\]
The computation of $\upi_k X(\infty\aO)$ then follows from
Corollary~\ref{cprops}.
    
    \item Equivalence in (2) follows from the cofibre sequence ${E\cF_{0}}_{+}\to S^0\to\EFt_0$.

    \item  Note that $\underline{f}_{C_2}$ be the $C_2$-Mackey functor with $\underline{f}_{C_2}(C_2/e)=\bF$ and $\underline{f}_{C_2}(C_2/C_2)=0$, i.e.\ the kernel of $\uF\to\ug_{C_2}$. Applying $\Phi^*_{H_0}$ to the short exact sequence of $C_2$-Mackey functors $0\to\underline{f}_{C_2}\to\uF\to\ug_{C_2}\to0$, and using Theorem \ref{tformula}, $\Phi^*_{H_0}\uF=\uF$ and $\Phi^*_{H_0}\ug_{C_2}=\ug_{(0)}$, gives
\[
0\to\Phi^*_{H_0}\underline{f}_{C_2}\to\uF\to\ug_{(0)}\to 0.
\] 

From~\cite{GK26}, we have 
\[
\upi_{k}(S^{n\aO}\wedge H\uF)\cong 
\begin{cases}
\ug_{(0)} & 0\leq k \leq n-1,\\
\uF & k=n,\\
0 & \text{otherwise},
\end{cases}
\]
where non-zero classes at $K_4/K_4$ level is given by $a^{n-k}_{\alpha_0}u^{k}_{\alpha_0}$ for $0\leq k \leq n$. Then the computation of $\upi_{k} F(\alpha_0;n)$ follows from the Mackey functor map $\upi_{k} X({n\aO}) \to \upi_{k}X(\infty \aO)$. 
\end{enumerate}\end{proof}

Having computed the homotopy Mackey functors of the fibre $F(\alpha_i; n)$, we now
use the cofibre sequence~\eqref{fibalpha_0} to compare the slices of
$X(n\alpha_i)$ and $X(\infty\alpha_i)$.

\begin{thm}\label{thm:Xinfslices}
  \begin{enumerate}
    \item $X(n\alpha_i)$ and $X(\infty\alpha_i)$ have the same slices in all levels
      $< n$. 
    \item The slices of $X(\infty\alpha_i)$ are concentrated in the even levels
      $2\ell$ with $\ell \geq 0$. 
      \item The slices of $X(\infty \alpha_i)$ are given by
 \[
P^{4m}_{4m}X(\infty \alpha_i)\cong
\begin{cases}
 H\ug_{(i)} & m=0,\\
S^{ \rho}\wedge  H\ug_{(i)}^{*} & m=1,\\
S^{m \rho}\wedge  H(\phi^*_{H_j}\uF^{*}\oplus\phi^*_{H_k}\uF^{*}\oplus\ug^{\oplus(2m-3)}) & m\ge2,
\end{cases}
\]
and 
\[
P^{4m+2}_{4m+2}X(\infty\alpha_i)\simeq
\begin{cases}
\Sigma H\ug_{(i)} & m=0,\\
S^{m\rho}\wedge\Sigma H(\phi^{\ast}_{H_j}\uf_{C_2}\oplus\phi^{\ast}_{H_k}\uf_{C_2}) & m\ge1,
\end{cases}
\]
where $i,j,k$ are distinct elements of the set $[2]$. 
  \end{enumerate}
\end{thm}

\begin{proof} Proof for $i=0$.

(1) Consider the  cofibre sequence \eqref{fibalpha_0}, by Proposition~\ref{prop:Xinf}, $F(\aO;n)$ is $n$-slice-connective. Hence, by \eqref{S:HY}, we have $F(\aO;n)\geq n$. Hence the desired equivalences follow from Lemma~\ref{lPequiv}. 

(2)  It is evident that $X(\infty\alpha_0)\geq 0$, and hence it has no negative slices. Consider the cofibre sequence
\[
\Sigma^{1+\alpha_0} X(\infty\alpha_0)\xrightarrow{\ u\ }X(\infty\alpha_0)\longrightarrow H\ug_{(0)},
\]
which tells us $\res_H X(\infty\alpha_0)$ has no odd slices for any proper subgroup $H \subset K_4.$

If $X(\infty \alpha_0)$ have odd slices, then it must be geometric.  Therefore, the slices occur in dimensions divisible by $|K_4|=4$, so they occur in even dimensions.

$(3)$ We start with analyzing $4m$-slices of $X(\infty \alpha_0).$ Note that $S^{-m\rho}\wedge P^{4m}_{4m}X(\infty \alpha_0)$ is a $0$-slice, by $\rho$-periodicity. Every $0$-slice is Eilenberg--Mac Lane space, $H\uM$ with
$\uM=\upi_0$. Hence we have
\[
P^{4m}_{4m}X(\infty \alpha_0)\simeq S^{m\rho}\wedge H \upi_{m\rho}X(\infty \alpha_0).
\]
Using Anderson duality, we have

\[
\upi_{m\rho}X(\infty \alpha_0)=\underline{H}^{K_4}_{m\rho}(\widetilde{E\{e, H_0\}};\uF)\cong
\mathrm{Hom}_{\uF}\bigl(\underline{H}^{(m+3)+(m-1)\bar\rho}_{K_4}(\widetilde{E\{e, H_0\}};\uF),\uF\bigr)
\]
as Mackey functors. Under this isomorphism, restriction for $\upi_{m\rho}X(\infty \alpha_0)$ is dual to transfer in cohomology, and transfer is
dual to restriction.

 The cofibre sequence $E\{e, H_0\}_+\to S^0\to\widetilde{E\{e, H_0\}}$ gives an exact sequence of Mackey functors and the coboundary map 
\[
\underline{H}^{(m+2)+(m-1)\bar\rho}_{K_4}(E\{e, H_0\}_+; \uF)\xrightarrow{\ \delta\ }\underline{H}^{(m+3)+(m-1)\bar\rho}_{K_4}(\widetilde{E\{e, H_0\}}; \uF)
\]
is an isomorphism of Mackey functors.

 By~\cite[Proposition~2.5]{GK26}, the cohomology group $\widetilde H^{(m+2)+(m-1)\bar\rho}_{K_4}
\bigl(E\{e,H_0\}_+;\uF\bigr)$ admits a basis of the following form:
\[
x_t=a_{\alpha_0}^{2m+1+t}u_{\alpha_0}^{-(m+2+t)}u_{\alpha_1-\alpha_2}^{\,m-1}a_{\alpha_2}^{2m-2-t}u_{\alpha_2}^{t},\qquad 0\le t\le2m-2 
\]
when $m\geq 1 $.
While at orbits $K_4/H_1$ and $K_4/H_2$, the group is $\bF$, generated by $a_\sigma^{4m-1}u_\sigma^{-2m-1}\in \widetilde H^{2m+1-\sigma(2m-2)}_{C_2}({EC_{2}}_{+};\uF)$ from \cite{BG21}. The restrictions are
\[
\res_{H_1}(x_t)=a_\sigma^{4m-1}u_\sigma^{-2m-1}\text{ and }
\res_{H_2}(x_t)=\delta_{t,2m-2}\,a_\sigma^{4m-1}u_\sigma^{-2m-1} \text{ for all }t.
\]
 They are equal for $m=1$.  Since
each generator is a restriction, Frobenius reciprocity implies
\[
\mathrm{tr}_{H_1}(\res_{H_1}x_0)=x_0\,\mathrm{tr}_{H_1}(1)=0,\qquad \mathrm{tr}_{H_2}(\res_{H_1}x_{2m-2})=x_{2m-2}\,\mathrm{tr}_{H_2}(1)=0 .
\]

Let $\uM_{m}:=\underline{H}^{(m+2)+(m-1)\bar\rho}_{K_4}(E\{e, H_0\}_+;\uF)$, so that $\upi_{m\rho}Y\cong\underline M^{*}_{m}\cong \Hom_{\mathbb{F}_2}(\uM_{m}, \mathbb{F}_2).$ Let us denote $\{x_0^*, \cdots x_{2m-2}^*\}$ be the basis of $\uM^*_{m}(K_4/K_4)$ dual to basis $\{x_0,x_1,\cdots,x_{2m-2}\}$. Similarly, $y^*$ is the basis of
$\uM^*_{m}(K_4/H_1)\cong\bF$ (respectively $\uM^*_{m}(K_4/H_2)$) dual to $y:= a_\sigma^{4m-1}u_\sigma^{-2m-1}$. Then
\[
\upi_{m\rho}Y(K_4/K_4)\cong\bF\langle x_0^*,\dots,x_{2m-2}^*\rangle\cong\bF^{2m-1},
\qquad
\upi_{m\rho}Y(K_4/H_1)\cong\upi_{m\rho}Y(K_4/H_2)\cong\bF\langle y^*\rangle .
\]
Dualization interchanges restrictions and transfers. Precisely, for $H\in\{H_1,H_2\}$,
\[
\res^{K_4}_H=(\mathrm{tr}^{K_4}_H)^* \text{ and } \mathrm{tr}^{K_4}_H=(\res^{K_4}_H)^*\qquad\text{in    }\upi_{m\rho}Y,
\]
where the maps on the right are the structure maps of $\underline M_m$. Hence the restrictions $\upi_{m\rho}Y(K_4/K_4)\to\upi_{m\rho}Y(K_4/H)$ vanish, because $\mathrm{tr}_{H_1}=\mathrm{tr}_{H_2}=0$ in $\uM_{m}$ and the transfers $\upi_{m\rho}Y(K_4/H)\to\upi_{m\rho}Y(K_4/K_4)$ are given by $y^*\mapsto y^*\circ\res_H$, that is,
\[
\mathrm{tr}^{K_4}_{H_1}(y^*)=x_0^*+x_1^*+\dots+x_{2m-2}^*,\qquad \mathrm{tr}^{K_4}_{H_2}(y^*)=x_{2m-2}^* .
\]
 For $m=1$,  we have $\ug_{(0)}^{*}$. For $m\ge2$, the two transfer images span a $2$-dimensional subspace giving
$\phi^*_{H_1}\uF^{*}\oplus\phi^*_{H_2}\uF^{*}$, and any complement is $\ug^{\oplus(2m-3)}$.

Thus we derive the $4m$-slice of $X(\infty \alpha_0)$ as

\[
P^{4m}_{4m}X(\infty \alpha_0)\cong
\begin{cases}
 H\ug_{(0)} & m=0,\\
S^{\rho}\wedge  H\ug_{(0)}^{*} & m=1,\\
S^{m \rho}\wedge  H(\phi^*_{H_1}\uF^{*}\oplus\phi^*_{H_1}\uF^{*}\oplus\ug^{\oplus(2m-3)}) & m\ge2.
\end{cases}
\]

Next we compute $4m+2$-slices of $X(\infty \alpha_0)$. For $P^{4m+2}_{4m+2}X(\infty \alpha_0)$, let us denote the spectrum $B_m:= S^{-m\rho}\wedge P^{4m+2}_{4m+2}Y$ which is a $2$-slice by \eqref{S:rho}.

By Slone's characterization of $2$-slices for group $K_4$ \cite[Thm.~5.9]{Slo22}, only $\upi_1B_m$ and $\upi_2B_m$ can be nonzero, all
restrictions of $\upi_2B_m$ are injective, and $\upi_1B_m(e)=0$.

Let us fix $Y=X(\infty \alpha_0).$ Since $\res_eY\simeq*$, we have $\upi_2B_m(K_4/e)=0$. Injectivity of the restrictions then forces $\upi_2B_m=0$,
first at $H_0$, $H_1$, $H_2$ and then at $K_4$. So $B_m\simeq\Sigma H\underline L_m$ with $\underline L_m=\upi_1B_m$.

Note that $\res_eY\simeq\res_{H_0}Y\simeq*$ and the slices of
$\res_{H_1}Y\simeq\widetilde {EC_2}\wedge H\uF$ are $\Sigma^mH_{C_2}\ug_{C_2}$ in level $2m$. Hence $\underline L_m$ is $0$ at $e$ and $H_0$ and $\bF$ at ${H_1}$ and ${H_2}$.

Since $\Sigma H\underline L_m$ is a $2$-slice, \cite[Prop.~5.7]{Slo22} shows that
$\underline L_m(K_4)\to\underline L_m(H_1)\oplus\underline L_m(H_2)=\bF^2$ is injective.

By the double coset formula, $\res^{K_4}_{H_1}\mathrm{tr}^{K_4}_{H_1}=2=0$ and
$\res^{K_4}_{H_1}\mathrm{tr}^{K_4}_{H_1}=\mathrm{tr}^{H_2}_e\res^{H_2}_e=0$. Together with injectivity this forces $\mathrm{tr}^{K_4}_{H_1}=0$, and likewise
$\mathrm{tr}^{K_4}_{H_2}=0$.

The automorphism of $K_4$ that swaps ${H_1}$ and ${H_2}$ fixes $\aO$, so it fixes $Y$. Hence the image of
$\underline L_m(K_4)$ in $\bF^2$ is $0$, the diagonal or $\bF^2$, which gives  $\underline{L}_m$ is one of
\[
\phi^{*}_{H_1}(\underline{f}_{C_2})\oplus \phi^{*}_{H_2}(\underline{f}_{C_2}),\qquad \ug_{(0)},\qquad \phi^*_{H_1}\uF\oplus\phi^*_{H_2}\uF .
\]

Consider the cofibre sequence $P_{4m+4}Y\to P_{4m+2}Y\to S^{m\rho}\wedge B_m$. Now apply $\upi_{m\rho+1}$ to get
\[
0 =\upi_{2}B_m \to \upi_{m\rho+1}P_{4m+4}Y\xrightarrow{\ f\ }\upi_{m\rho+1}Y\to \upi_1B_m\to 0.
\]

Restriction maps at level $K_4/K_4$ to the sum of its values at $K_4/{H_1}$ and $K_4/{H_2}$ give
\[
\xymatrix@C=1.6pc{
0\ar[r]&\upi_{m\rho+1}P_{4m+4}Y(K_4/K_4)\ar[d]\ar[r]^-{f}&\upi_{m\rho+1}Y(K_4/K_4)\ar[d]^{(\res_{H_1},\res_{H_2})}\ar[r]^-{g}&\upi_1B_m(K_4/K_4)\ar[d]^{(\res_{H_1},\res_{H_2})}\ar[r]&0\\
0\ar[r]&0\ar[r]&\bF^2\ar[r]^-{\cong}&\bF^2\ar[r]&0 .}
\]

The left vertical map has target $0$, so by commutativity of the left square
the restriction of any element of $\im(f)$ is zero.

 By Slone \cite[Prop.~5.7]{Slo22}, the right vertical map is injective. Let
$x\in\upi_{m\rho+1}Y(K_4/K_4)$ with $(\res_{H_1},\res_{H_2})(x)=0$. By commutativity of the right square, $(\res_{H_1},\res_{H_2})(g(x))=0$. Injectivity
gives $g(x)=0$, and exactness gives $x\in\im(f)$.

Hence $\im(f)=\ker(\res_{H_1},\res_{H_2})$. Since $g$ is surjective, it follows that
\[
\upi_1B_m(K_4/K_4)\;\cong\;\upi_{m\rho+1}Y(K_4/K_4)\big/\ker(\res_{H_1},\res_{H_2})\;\cong\;
\im\Bigl(\upi_{m\rho+1}Y(K_4/K_4)\xrightarrow{(\res_{H_1},\res_{H_2})}\bF^2\Bigr).
\]

Consider the cofibre sequence induced by multiplication by \(u_{\alpha_0}\).
\[
\Sigma Y \to Y \to H\ug_{(0)}.\]
 From the associated long exact sequence, we derive that $\upi_{m\rho}Y \to \upi_{m\rho+1}Y$ is multiplication by $t=u_{\alpha_0}a_{\alpha_0}^{-1}\in\pi^{K_4}_1Y$.
 and it is an isomorphism of Mackey functors. Therefore, we have the commutative diagram:
 \[
 \xymatrix{\upi_{m\rho}Y(K_4/K_4) \ar[d]^{\res^{K_4}_H}\ar[r]^t & \upi_{m\rho+1}Y(K_4/K_4)\ar[d]^{\res^{K_4}_H} \\ \upi_{m\rho}Y(K_4/H)  \ar[r] & \upi_{m\rho+1}Y(K_4/H) }
 \]
 for $H \in \{{H_1}, {H_2}\}.$ For $m\ge 1,$ the restriction maps $\res_{H_1}, \res_{H_2}$ for the Mackey functor $\upi_{m\rho}Y$ are zero. Hence, $\upi_1 B_m(K_4/K_4)=0.$  Therefore, we obtain
\[
P^{4m+2}_{4m+2}X(\infty\aO)\simeq
\begin{cases}
\Sigma H\ug_{(0)} & m=0,\\
S^{m\rho}\wedge\Sigma H(\phi^{*}_{H_1}\underline{f}_{C_2}\oplus \phi^{*}_{H_2}\underline{f}_{C_2})) & m\ge1.
\end{cases}\qedhere
\]\qedhere
\end{proof}
We now describe the \(n\)-slices of \(X(n\alpha_0)\). By the symmetry of \(K_4\), which permutes the pairs \((H_i,\alpha_i)\), the corresponding statements for \(X(n\alpha_i)\) are obtained by replacing \(H_0\) with \(H_i\). Hence, we restrict our attention to the case \(i=0\).
\begin{example}
The \(1\)-slice of \(X(\alpha_0)\) is \(\Sigma H\uF\), while the \(2\)-slice of \(X(2\alpha_0)\) is \(\Sigma^{1+\alpha_0}H\uF\). 
\end{example}

The cases \(n\geq 3\) are discussed in the following theorem.

\begin{thm}\label{pnslice}
For $n\ge 3$, we have
\begin{myeq}{\label{nslice}}
P^{n}_{n} (X(n\alpha_0))\simeq \begin{cases}
    \Sigma^{\frac{n}{4} \rho} H(\underline{N}_{H_0}\oplus \ug^{\frac{n}{2}-2}) & n\equiv 0 \pmod 4,\\
  \Sigma^{\frac{n-1}{4}\rho+1}H\underline{K}_{H_0} & n\equiv 1 \pmod 4,\\
    \Sigma^{\frac{n-2}{4}\rho}X_0 & n\equiv 2 \pmod 4,\\
\Sigma^{\frac{n+1}{4}\rho-1}H\underline{D}_{H_0} & n\equiv 3 \pmod 4.  
\end{cases}    
\end{myeq}
where the Mackey functors  $\underline{N}_{H_0}$, $\underline{K}_{H_0}$ and $\underline{D}_{H_0}$ are given as follow.
\begin{myeq}{\label{makK}}
\begin{array}{|c|c|c|}
\hline
\xymatrix@C=0.9cm@R=1.1cm{
  & \F \ar@/_.5pc/[dl]^{1} \\
  \F
  & \F\ar@/_.5pc/[d]^{1} \ar@/_.5pc/[u]^{1}
  & \F\ar@/^.5pc/[dl]^{1}
  \ar@/_.5pc/[ul]^{1}\\
  & \F\ar@/^.5pc/[ul]^{1}
}
&
\xymatrix@C=0.9cm@R=1.1cm{
 & 0 & \\
0 & \F\ar@/^.5pc/[d] & \F\ar@/^.5pc/[dl]\\
 & \F &
}
&
\xymatrix@C=0.9cm@R=1.1cm{
 & 0 & \\
\F & 0 & 0 \\
 & \F\ar@/^.5pc/[ul] &
}
\\[0.5cm]
\hline
\underline{N}_{H_0}
&
\underline{K}_{H_0}
&
\underline{D}_{H_0}
\\
\hline
\end{array}
\end{myeq}

The spectrum $X_0:=P^{2}_{2}X((n\alpha_0)-\frac{n-2}{4}\rho)$ is $2$-slice whose  Mackey functor valued homotopy groups are concentrated in degree $1$ and $2$, $\upi_{1} X_0\cong \phi^{\ast}_{H_1}\uf_{C_2} \oplus \phi^{\ast}_{H_2}\uf_{C_2}$ and $\upi_{2} X_0\cong \underline K_{H_0}$. Hence there is a cofibre sequence
\[
\Sigma^{2+\frac{n-2}{4}\rho} H\underline{K}_{H_0} \to P^{n}_{n}X(n\aO)\to \Sigma^{1+\frac{n-2}{4}\rho}H(\phi^{\ast}_{H_1}\uf_{C_2} \oplus \phi^{\ast}_{H_2}\uf_{C_2}).
\]
Moreover, apart from the \(n\)-slice, there are no odd-dimensional slices.
\end{thm}

\begin{proof} In the proof, we consider $n$ is of the form \(n=4\ell+t\), for $t=0,1,2,3$ and $\ell\geq 0$.

We begin the proof that, apart from $n$, there are no odd slices. Let $m\neq n$ be an odd integer, then  $\res_H P^{m}_{m}X(n\alpha_0)\simeq P^{m}_{m}\res_H X(n\alpha_0)\simeq \ast$ for any $2$-order subgroup $H$ of $K_4$. Therefore, the slice seen only by $K_4$ and hence lives in a level divisible by $4$. So the only possible odd slice is in level $n$, which occurs when $n$ is odd.

Now we consider the computation of the $n$-slice for $X(n\aO)$. We analyze this case-by-case. 

{\bf Case I} Let $n=4\ell$. From~\eqref{S:01}, the $0$-slice is of the form $H\uM$, where $\uM\simeq\upi_{0}(X(4\ell\alpha_0-\ell\rho))$. The value of the Mackey functor \(\uM\) at the orbit \(K_4/K_4\) from~\cite{GK26} is given by
\[
\uM(K_4/K_4)\cong
\F\left\{
a^{s}_{\alpha_0}u^{2\ell-s}_{\alpha_0}v^{\ell}_{1}
\Sigma^{-1}
\frac{1}{a^{1+s}_{\alpha_2}u^{2\ell-1-s}_{\alpha_2}}
\;\middle|\;
0\leq s\leq 2\ell-2
\right\}.
\]
Note that the classes with \(s\neq 0\) generate copies of \(\ug\), whereas the class with \(s=0\) generates \(\underline{N}_{H_0}\).

{\bf Case II.}  For $n= 4\ell+1$, we have
  \[
  P^{4\ell+1}_{4\ell+1}X((4\ell+1)\alpha_0) \simeq \Sigma^{\ell\rho}P^{1}_{1}X((4\ell+1)\alpha_0-\ell\rho)
  \]
By~\eqref{S:01}, we have $P^{1}_{1}\bigl(X((4\ell+1)\alpha_0-\ell\rho)\bigr)
\simeq \Sigma H\uM,$ where \(\uM\) is a Mackey functor with injective restriction maps. To determine \(\uM\), we restrict the above slice to the order-two subgroups and apply Proposition~\ref{C_2slices}.

\begin{align*}
\res _{H_t}P^{1}_{1}X((4\ell+1)\alpha_0-\ell\rho)  \simeq\begin{cases}
    \Sigma^{1}H_{C_2}\uf_{C_2}, & t=0,\\
 \Sigma^{1}H_{C_2}\uF, & t\neq 0.
\end{cases} 
\end{align*}

These computations force $\uM\cong \underline{K}_{H_0}$.

{\bf Case III.}~ For $n=4\ell+2$, 
the spectrum $Y= X((4\ell+2)\alpha_0-\ell\rho)\in [-4\ell, 4\ell]$ from Corollary~\ref{bound0} and let $X_0= P^2_2Y$ denote $2$-slice of the respective spectrum. Note that $Y$ has no odd slice.  By \eqref{S:01}, we have $P^{4\ell+4}_{4\ell+4}X(n\aO) \simeq S^{(\ell+1)\rho}\wedge H\uM$, where $\uM\cong \upi_{(\ell+1)\rho}X(n\aO).$ Since \(\uM(K_4/H_1)=\uM(K_4/H_2)=\uM(K_4/e)=0\), we obtain $
P_4^4Y\simeq S^\rho\wedge H\uM
\simeq S^{1+\alpha_0}\wedge H\uM.$ Moreover, we have the following exact sequence:
\[
\upi_2X_0\to\upi_{1}P^4_{4}Y \to\upi_{1}Y \xrightarrow{\ g\ }\upi_1X_0\to0,
\]
which is induced from the sequence of isomorphisms along with the cofibre sequence $P_4 \to P_2\to X_0$.

To compute $\upi_{1}P^4_{4}Y$, let us consider the cofibre sequence
\[
{K_4/H_0}_+\to S^0\to S^{\aO}
\]
which gives $\pi_{1}P^4_{4}Y\cong \pi_{-{\aO}}(H\uM)= \coker(\tr^{K_4}_{H_0}(\uM))$.

For $k\ge1$, $\pi^{K_4}_{k\rho}X(n\aO)\cong\widetilde{H}^{\,n\aO-k\rho}_{K_4}(S^0;\uF)$ has basis
\[
\{a_{\aO}^{\,n-2k-1-s}\,k_1^{\,k}u_{\aO}^{\,s+1}\Sigma^{-1}\frac{1}{a_{\aII}^{2k-s}u_{\aII}^{s}}:1\le s\le\min(2k-1,n-2k-1)\}
\] from~\cite{GK26}. 

Hence for $n=4\ell+2$ and $k=\ell+1,$ we have $\uM(K_4/K_4)\cong \F^{2\ell-1}$ and transfer $\tr^{K_4}_{H_0}(\uM)$ is zero. This gives $\pi_{1}P^4_{4}Y= \uM(K_4/K_4)\cong \F^{2\ell-1}.$
 Using the computation in \cite{HK17} of the Poincare series of the $RO(K_4)$-Bredon cohomology of $S^0$, it turns out that $\pi_1(Y)\cong \F^{2\ell-1}.$ This forces $\upi_1(X_0)(K_4/K_4)=0$. Hence $\upi_{1}(X_0)\cong \phi^{\ast}_{H_1}\uf_{C_2} \oplus \phi^{\ast}_{H_2}\uf_{C_2}.$  While the computation $\upi_{2}X_{0}=\underline{K}_{H_0}$ follows from the restriction to the order $2$-subgroup along with \cite[Theorem~5.9]{Slo22}.

{\bf Case IV.}~ For $n=4\ell-1$, the result follows by a similar argument to that in $n\cong 1 \pmod4$, together with~\eqref{S:01} for the $-1$-slice.
\end{proof}





\begin{example} Consider the spectrum $X(5\alpha_0)=S^{5\alpha_0}\wedge H\uF$. Using Corollary~\ref{bound0}, the nontrivial slices of $X(5\alpha_0)$ lies in the interval $[0,6]$. From Theorem~\ref{thm:Xinfslices} and Theorem~\ref{pnslice}, we can obtain all the slices up to level $5$.

For the computation of the remaining slice, $P^{6}_{6}(X(5\alpha_0))\simeq \Sigma^{\rho}P^{2}_{2}(X(5\alpha_0-\rho))$. So it is sufficient to compute the $2$-slice of $X(5\alpha_0-\rho)$. 

First, consider the restriction of the $2$-slice to the order $2$-subgroup and apply Proposition~\ref{C_2slices}; we have

\[
P^{2}_{2}(\res_{H_i}(X(5\alpha_0-\rho))) \simeq \begin{cases}
    P^{2}_{2}(\Sigma^{3-2\sigma}H_{C_2}\uF)\simeq\Sigma H_{C_2}\ug_{C_2} & i=0\\
    P^{2}_{2}(\Sigma^{3\sigma-2}H_{C_2}\uF)\simeq \ast & i=1,2.
\end{cases}
\]

Hence, $\upi_{2}(P^{2}_{2}X(5\alpha_0-\rho))$ is zero at orbit $K_4/H_i$ for $i=0,1,2$. Hence by \cite[Theorem~5.9]{Slo22}, $\upi_{2}(P^{2}_{2}X(5\alpha_0-\rho))=0$ at orbit $K_4/K_4$ orbit. 

To compute $\pi_{1}P^{2}_{2}$ at level $K_4/K_4$,  consider the following long exact sequence
\begin{myeq}{\label{long1}}
\xymatrix@C=0.8pc{
\cdots\ar[r]&\pi_{1}P^{2}_{2}X(5\alpha_0-\rho)\ar[r]^-{f}&\pi_{1}P^2X(5\alpha_0-\rho)\cong \F\{u^{3}_{\alpha_0}v_2\Sigma^{-1}\frac{1}{a_{\alpha_1}u_{\alpha_1}}\}\ar[r]&\pi_{1}P^{1}X(5\alpha_0-\rho)\ar[r] &\cdots}
\end{myeq}

We claim that $\pi_{1}P^{1}X(5\alpha_0-\rho)=0$. To see this, consider the following long exact sequence
\[
\cdots\to \pi_{1}P^{k}_{k}X(5\alpha_0-\rho) \to \pi_{1}P^{k}X(5\alpha_0-\rho) \to \pi_{1}P^{k-1}X(5\alpha_0-\rho) \cdots
\]
for $-4\le k\le 1$.

Using $P^{k}_{k}X(5\alpha_0-\rho)\simeq \Sigma^{-\rho}P^{k+4}X(5\alpha_0)$, for  $-4\le k\le 1$ and Theorem~\ref{thm:Xinfslices}, we have $\pi_{1} P^{a}X(5\alpha_0-\rho)\cong \pi_{1} P^{b}X(5\alpha_0-\rho)$ for $-5\le a,b\le 1$. But $P^{-5}X(5\alpha_0-\rho)\simeq \ast$. Hence $\pi_{1}P^{1}(5\alpha_0-\rho)=0$.

In \eqref{long1}, $f$ is a nonzero map, hence $\pi_{1}P^{2}_{2}X(5\alpha_0-\rho)\cong \F$. So the $2$-slice is
\[
P^{2}_{2}X(5\alpha_0-\rho) \simeq \Sigma {H{\phi}^{\ast}_{H_0}}\uF.
\]
Thus, we have the following slice tower for $X(5\alpha_0)$: 
\[
\xymatrix{
 P^{6}_{6}\simeq\Sigma^{\rho+1}H\phi^{\ast}_{H_0}\uF  \ar[r] &  P^6 X(5\alpha_0)\simeq X(5\alpha_0) \ar[d] \\
 P^{5}_{5}\simeq
\Sigma^{\rho+\aO} H\uf \ar[r] &  P^5 X(5\alpha_0) \ar[d] \\
P^{4}_{4}\simeq\Sigma H\ug \ar[r] &  P^4 X(5\alpha_0) \ar[d] \\ 
P^{2}_{2}\simeq\Sigma H \ug_{(0)} \ar[r]& P^2X(5\alpha_0)\ar[d] \\ & P^0 X (5\alpha_0)\simeq P^0_0 X(5\alpha_0)\simeq H\ug_{(0)}.}
\]

\end{example}

\begin{example} Next consider the spectrum $X(10\alpha_0)=\Sigma^{10\alpha_0}H\uF\in [0,16]$ and it has the following slice tower. 
\[
\xymatrix{
P^{16}_{16}\simeq\Sigma^{4\rho}H\phi^{\ast}_{H_0}\uF \ar[r] &  P^{16}  X(10\alpha_0)\simeq X(10\alpha_0) \ar[d] \\
P^{14}_{14}\simeq \Sigma^{3\rho+1}H{\phi^{\ast}_{H_0}\uf_{C_2}}  \ar[r] &  P^{14}  X(10\alpha_0) \ar[d] \\
P^{12}_{12} \simeq \Sigma^{3\rho}H\ug^{\oplus3} \ar[r] &  P^{12}  X(10\alpha_0) \ar[d] \\
P^{10}_{10}  \ar[r] &  P^{10}  X(10\alpha_0) \ar[d] \\
  P^8_8\simeq S^{2\rho} \wedge H(\phi^{\ast}_{H_1}\uF^{\ast} \oplus \phi^{\ast}_{H_2}\uF^{\ast}\oplus \ug) \ar[r] &  P^8  X(10\alpha_0) \ar[d] \\
 P^6_6 \simeq S^\rho\wedge \Sigma H(\phi^{\ast}_{H_1}\uf_{C_2}\oplus\phi^{\ast}_{H_2}\uf_{C_2}) \ar[r] &  P^6 X(10\alpha_0) \ar[d] \\
P^4_4\simeq S^{\rho}\wedge H\ug^{\ast}_{(0)}  \ar[r] &  P^4 X(10\alpha_0) \ar[d] \\  
P^2_2\simeq \Sigma H\ug_{(0)} \ar[r]& P^2X(10\alpha_0)\ar[d] \\ & P^0 X (10\alpha_0)\simeq P^0_0 X(10\alpha_0)\simeq H\ug_{(0)}.}
\] 
 The slice at level $10$ is given by the following extension.
 \[
 \Sigma^{2+2 \rho} H\underline K_{H_0} \to P^{10}_{10} X(10\alpha_0) \to \Sigma^{1+2\rho} H(\phi^{\ast}_{H_1}\uf_{C_2}\oplus\phi^{\ast}_{H_2}\uf_{C_2}).
 \]
 
\end{example}

 \section{Defect and spherical slice}
In this section, we begin by introducing the necessary definitions and terminology and then use these preliminaries to establish our main results, Theorem~\ref{main} and Theorem~\ref{mainodd}, concerning spherical slices. Later we reprove the result from~\cite{GY20} about $n$-slices of the spectra of the form $\Sigma^{n}H\uF$ for $n\geq 0$. Recall that, for any \(V\in RO(K_4)\), we say that \(\Sigma^{V}H\uF\) has a \emph{spherical slice} if there exists some \(U\in RO(K_4)\) such that
\[
P^{\dim V}_{\dim V}\Sigma^{V}H\uF\simeq \Sigma^{U}H\uF.
\]

 Throughout this section, we assume $V\in RO(K_4)$ have dimension $n$.

 For integers $l\le u$, let $\pr{x}{\ell}{u}=\max(l,\min(u,x))$ be the nearest point of $[l,u]$.
\begin{defn}\label{ddefect}
For a representation \(V\) of dimension \(n\), define the \emph{target dimensions} and the \emph{defect vector} as follows.
\[
\delta_{H_i}(V)=\pr{d_{H_i}(V)}{\lceil n/2\rceil}{\lfloor n/2\rfloor+2},\;
\Delta(V)=(\Delta_0(V),\Delta_1(V),\Delta_2(V)), \;\Delta_i(V)=\delta_{H_i}(V)-d_{H_i}(V) 
\]
for $i=0,1,2$ and $d_{H_i}(V)=\dim V^{H_i}$.

For Simplicity, we write $\sum_ {\Delta(V)}:=\Delta_0(V)+\Delta_1(V)+\Delta_2(V)$ and $[2]$ to denote the set $\{0,1,2\}$. If $\sum_ {\Delta(V)}$ is even. Then for each $i\in [2]$, define
\[
x_i(V):=\tfrac12(\Delta_j(V)+\Delta_k(V)-\Delta_i(V))
\]
for $i,j,k$ be distinct element of $[2]$.
\end{defn}

\begin{remark}
By the above definition, \(x_i(V)\)'s are integers if and only if \(\sum_{\Delta(V)}\) is an even integer.
\end{remark}
The following proposition describes the behaviour of slices upon restriction to subgroups of order two.

\begin{prop}\label{ptarget} For $V\in RO(K_4)$ with $\dim V=n$, then we have  \[
\res_{H_i}P^n_nX(V)\simeq S^{\delta_{H_i}+(n-\delta_{H_i})\sigma}\wedge H_{C_2}\uF
\]
for $i=0,1,2$.
\end{prop}
\begin{proof}
Since restriction commutes with the slice tower, so $\res_{H_i}P^n_nX(V)$ is $n$-slice of the $C_2$-spectrum
\[
\res_{H_i}X(V)=S^{d_{H_i}(V)+(n-d_{H_i}(V))\sigma}\wedge H\uF.
\] Over $C_2$, this slice is the representation sphere whose fixed-point dimension is the projection of $d_{H_i}(V)$ onto the interval $[\lceil n/2\rceil,\lfloor n/2\rfloor+2]$. 
\end{proof}

Since our focus is on the $n$-slice, we next determine the conditions under which the upper or lower bound, expressed in terms of $n(V)$ and $k(V)$, respectively, is equal to $n$; see~\eqref{boundn} and~\eqref{boundk}.

\begin{lemma}
Let $G$ be a finite group and $\uM$ a $G$-Mackey functor with $\uM(G/e)\neq0$.
Then for $V\in RO(G)$ with  $\dim V=n$, we have $n(V)\le n\le k(V)$. 

Moreover $n(V)<n$ if and only if $|H|(d_H(V)+b_H(\uM))<n$ for some subgroup $H\neq e$ and $k(V)>n$ if and only if
$|H|(d_H(V)-c_H(\uM))>n$ for some subgroup $H\neq e$.
\end{lemma}
\begin{proof}
If $\uM(G/e)\ne0$ then the underlying spectra of $H\uM$ and $H\uM^\vee$ are nonzero Eilenberg--MacLane spectra in degree $0$,
which implies $b_e=c_e=0$. As $d_e=n$, the $e$-terms in Theorem~\ref{trange} both equal $n$, giving the inequalities. The strict versions follow by inspecting the remaining terms. Alternatively, $\res^G_e$ preserves $\ge m$ and $\le m$,
and $\res^G_e X(V)\simeq\Sigma^nH\uM(G/e)$ is a nonzero $n$-slice, so any $m$ with $X(V)\ge m$ (resp.\ $\le m$)
satisfies $m\le n$ (resp.\ $m\ge n$).
\end{proof}

In Definition~\ref{ddefect}, we introduce the defect vector \(\Delta(V)\), whose \(i\)-component is the difference between the target dimension \(\delta_{H_i}(V)\) and the dimension of the \({H_i}\)-fixed point subspace of \(V\). The following lemma establishes that the sign on $\Delta(V)$ (on each component \(\Delta_{i}(V)\)) is determined by whether \(X(V)\leq n\) or \(X(V)\geq n\).

 \begin{lemma}\label{lnonnegdef} For \(V\in RO(K_4)\), the defect vector \(\Delta(V)\) satisfies the following properties:
\begin{enumerate}
    \item If \(X(V)\leq n\), then $\Delta_i(V)= \max(0, \lceil n/2\rceil -d_{H_i}(V))$ for $i \in [2]$, hence \(\Delta(V)\geq 0\).
    \item If \(X(V)\geq n\), then $\Delta_{i}(V)= \min(0,\,\lfloor n/2\rfloor+2-d_{H_i}(V))$ for $i \in \{0,1, 2\}$, hence \(\Delta(V)\leq 0\).
\end{enumerate}
 \end{lemma}
 \begin{proof} We prove (1), as the proof of (2) follows analogously.
 
     By Corollary~\ref{pk4range}, under the
hypothesis $X(V)\le n$, we have $n\geq k(V)$, in particular, \ $d_{H_i}(V)\le n/2+2$. As $d_{H_i}(V)\in\mathbb Z$, this means
$d_{H_i}(V)\le\lfloor n/2\rfloor+2$. Hence $\min(\lfloor n/2\rfloor+2,d_{H_i}(V))=d_{H_i}(V)$, which gives $\delta_{H_i}(V)=\max\big(\lceil n/2 \rceil, d_{H_i}(V))$. Therefore, $\Delta_{i}(V) \ge 0$ for $i=0,1,2.$ 
 \end{proof}

In Section~4, we discussed the maps induced by multiplication by the classes \(u_{\alpha_i}\) and \(v_i\), for \(i=0,1,2\). For an arbitrary \(V\in RO(K_4)\), these maps give rise to four possible moves.

\begin{prop}\label{pmoves}
For \(V\in RO(K_4)\), let \(i\in[2]\) and set \(X(V)=S^V\wedge H\uF\). Smashing the classes
\[
u_{\alpha_i}\colon\Sigma^{1-\alpha_i}H\uF\to H\uF
\qquad\text{and}\qquad
v_i\colon\Sigma^{\tau_i}H\uF\to H\uF
\]
with the representation sphere \(S^V\) yields the following maps, together with their respective cofibres and fibres.
\begin{enumerate}[label=(M\arabic*)]
\item \label{raiseu} $X(V+1-\alpha_i)\xrightarrow{\ u_{\alpha_i}\ }X(V)$, with cofibre $\Sigma^VH\ug_{(i)}$.
\item \label{loweru} $X(V)\xrightarrow{\ u_{\alpha_i}\ }X(V+\alpha_i-1)$, with fibre $\Sigma^{V-2}H\ug_{(i)}$.
\item \label{raisev}$X(V+\tau_i)\xrightarrow{\ v_i\ }X(V)$, with cofibre $\Sigma^VC(v_i)$.
\item  \label{lowerv}$X(V)\xrightarrow{\ v_i\ }X(V-\tau_i)$, with fibre $\Sigma^{V-\tau_i-1}C(v_i)$.
\end{enumerate}
All four maps preserve $\dim$. The  moves \ref{raiseu}, \ref{loweru} change $(d_{H_j}(V),d_{H_k}(V),d_{K_4}(V))$ by $\pm(1,1,1)$ and fix $d_{H_i}(V)$.
The moves \ref{raisev}, \ref{lowerv} change $(d_{H_i}(V),d_{K_4}(V))$ by $\pm(2,1)$ and fix $d_{H_j}(V),d_{H_k}(V)$.
\end{prop}

\begin{proof}
By definition, $u_{\alpha_i}$ and $v_i$ have cofibres $H\ug_{(i)}$ and $C(v_i)$, and $S^W\wedge(-)$ preserves cofibre sequences.
\end{proof}
Additionally, bounds on these fibres and cofibres are established in Proposition~\ref{psusgi} and Proposition~\ref{psuscvi}.

Finally, using the terminology introduced above, we characterize the spherical slice of \(X(V)\).

\begin{thm}{\label{main}}
\begin{enumerate}
    \item  Suppose $X(V)\le n$.  If $\sum_{\Delta(V)}$ is even, and $x_0(V), x_1(V), x_2(V)\ge 0$, there exists unique \[
U=V+x_0(V)(1-\aO)+x_1(V)(1-\aI)+x_2(V)(1-\aII)
\]
such that $P^n_nX(V)\simeq X(U).$
\item Suppose $X(V)\ge n$.  If $\sum_{\Delta(V)}$ is even, and $x_0(V), x_1(V), x_2(V)\le 0$, there exists unique \[
U=V+x_0(V)(1-\aO)+x_1(V)(1-\aI)+x_2(V)(1-\aII)
\]
such that $P^n_nX(V)\simeq X(U)$.
\end{enumerate}

\end{thm}

\begin{proof}
    (1) \textbf{Step 1.}
    Consider the representation \(U\in RO(K_4)\) constructed from \(V\in RO(K_4)\) as follows:
\[
U=V+x_0(V)(1-\aO)+x_1(V)(1-\aI)+x_2(V)(1-\aII).
\]
Then, for \(\{i,j,k\}=[2]\), we have
\begin{myeq}{\label{eqdi}}
d_{H_i}(U)=d_{H_i}(V)+x_j(V)+x_k(V)=\delta_{H_i}(V).   
\end{myeq}

By the definition of \(\delta_{H_i}(V)\), we have
\[
\left\lceil\frac n2\right\rceil
\le \delta_{H_i}
\le \left\lfloor\frac n2\right\rfloor+2.
\]
It follows that \(2\delta_{H_i}(V)\ge n\) and \(2\delta_{H_i}(V)-4\le n\). Consequently, by \eqref{eqdi}, we obtain
\[
0\le e_i(U)=2\delta_{H_i}(V)-n\le 4.
\]
Hence, by Corollary~\ref{pk4range}, \(X(U)\) is an \(n\)-slice.

\textbf{Step 2. }
Suppose that \(x_i(V)\ge 0\) for all \(i=0,1,2\). Then repeated application of the maps induced by \(u_{\alpha_0}\), \(u_{\alpha_1}\), and \(u_{\alpha_2}\) gives a map \(X(U)\to X(V)\). 

Let \(V_s\) be an intermediate representation, and let $\Delta_i(V_s)=\delta_{H_i}(V)-d_{H_i}(V_s)$
denote its remaining defect. Then
\[
\Delta_i(V_s)=x_j(V_s)+x_k(V_s),
\]
where \(x_j(V_s),x_k(V_s)\ge0\) count the remaining raises along \(\alpha_j\) and \(\alpha_k\), respectively. Suppose that \(x_i(V_s)>0\). We consider the map
\[
u_{\alpha_i}\colon X(V_s+1-\alpha_i)\longrightarrow X(V_s),
\]
whose cofibre is \(C(V_s+1-\alpha_i)\simeq \Sigma^{V_s}H\ug_{(i)}\).

For \(j\ne i\), we have
\[
\Delta_j(V_s)=x_i(V_s)+x_k(V_s)>0.
\]
Thus, $d_{H_j}(V)\le d_{H_j}(V_s)<\delta_{H_j}(V).$
It follows that \(\Delta_j(V)>0\), and hence \(\delta_{H_j}(V)=\lceil n/2\rceil\) by Lemma~\ref{lnonnegdef}. Therefore,
\begin{myeq}{\label{ej}}
2d_{H_j}(V_s)
\le 2\left\lceil\frac n2\right\rceil-2
\le n-1.
\end{myeq}
The same argument applies to \(H_k\).

Recall that \(e_i=2d_{H_i}(V)-n\), so that $4d_{K_4}(V)=n+\sum_{i=0}^{2}e_i.$

By the preceding argument~\eqref{ej}, $e_j(V_s),e_k(V_s)\le -1.$ Moreover, since \(d_{H_i}(V_s)\le\delta_{H_i}(V)\le \lfloor n/2\rfloor+2\), we have \(e_i(V_s)\le4\). Consequently, $4d_{K_4}(V_s)\le n+2,$ and therefore
\[
4\bigl(d_{K_4}(V_s)-1\bigr)\le n-2\le n-1.
\]
Combining these inequalities with the criterion for \(\Sigma^{V_s}H\ug_{(i)}\le m\) from Proposition~\ref{psusgi}, we obtain
\[
\Sigma^{V_s}H\ug_{(i)}\le n-1.
\]
Hence,
\[
P^n_nX(V_s+1-\alpha_i)\simeq P^n_nX(V_s).
\]

We now iterate this procedure. Start with \(V_0=V\) and \(x(V_0)=(x_0(V),x_1(V),x_2(V))\ge 0\). Whenever \(x_i(V_s)>0\), define
\[
V_{s+1}=V_s+1-\alpha_i.
\]
This raise leaves \(d_{H_i}\) unchanged and increases \(d_{H_j}\) and \(d_{H_k}\) by \(1\). Thus, in the order \((H_i,H_j,H_k)\),
\[
\Delta(V_{s+1})
=\Delta(V_{s})-(0,1,1),
\qquad
x(V_{s+1})=x(V_{s})-\mathbf e_i\ge 0,
\]
where \(\mathbf e_i\) denotes the \(i\)-th standard basis vector. Since \(x_0(V_{s})+x_1(V_{s})+x_2(V_{s})\) decreases by \(1\) at each step, after $x_0(V)+x_1(V)+x_2(V)$ steps, we reach \(x(V_{s_0})=0\). Consequently, \(\Delta{(s_0)}(V_{s_0})= 0\), and hence \(V_{s_0}=U\). Since each step induces an equivalence on the \(n\)-th slice, we conclude that
\[
P^n_nX(V)\simeq X(U).
\] 

(2) Proof of (2) is analogous to (1). In this case, we get a map $X(V) \to X(U)$ (map~\ref{loweru}) by repeatedly applying $u_{\aO}, u_{\aI}, u_{\aII} $, and we compute the fibres inductively. 
\end{proof}

When $\sum_{\Delta(V)}$ is odd, the following results hold.
\begin{thm}\label{mainodd} For $V\in RO(K_4)$, with $\dim V=n$ and $\sum_{\Delta(V)}$ is odd.
 
 (a) Suppose $X(V)\le n$ and
$x_{0}(V),x_{1}(V),x_{2}(V)\ge0$. Let us assume that $\Delta_i(V)$ is odd for some $i\in[2]$. Then $P^n_nX(V)$ is not spherical, and there is a cofibre sequence
\begin{myeq}{\label{odd1}}
X(\widehat{V'})\xrightarrow{\ v_i\ }P^n_nX(V)\to P_n(\Sigma^{\widehat{V'}-\tau_i}C(v_i)),
\end{myeq}
in which $X(\widehat{V'})$ is an $n$-slice. Here $ \widehat{V'}= V+\sum_{\ell=0}^2t_\ell\tau_\ell+\min(\Delta_j, \Delta_k)(1-\alpha_i)$ for some $t_\ell$, and $j$ and $k$ are the two distinct elements of $[2]\setminus\{i\}$.

(b) Suppose $X(V)\ge n$ and $x_{0}(V),x_{1}(V),x_{2}(V)\le0$. Assume $\Delta_i(V)$ is odd. Then $P^n_nX(V)$ is not spherical, and there is a cofibre sequence
\begin{myeq}{\label{odd2}}
P^n(\Sigma^{\widehat{V'}-1}C(v_i))\to  P^n_nX(V) \xrightarrow{\ v_i\ } X(\widehat{V'}),
\end{myeq}
in which $X(\widehat{V'})$ is an $n$-slice. Here $ \widehat{V'}= V-\sum_{\ell=0}^2t_\ell\tau_\ell-\min(|\Delta_j|, |\Delta_k|)(1-\alpha_i)$ for some $t_\ell,$ and $j$ and $k$ are the two distinct elements of $[2]\setminus\{i\}$.
\end{thm}
\begin{proof}We prove the result for $i=0$; the other cases follow similarly.

  Since \(\Delta_0(V)\) is odd, we have \(\Delta_1(V)\equiv\Delta_2(V)\pmod{2}\). Recall that $\tau_0=1+\aO-\aI-\aII$. Therefore  $\Delta_1(V+\tau_0)$ and $\Delta_2(V+\tau_0)$ are unchanged. But $\Delta_0(V+\tau_0)$,
\[
\Delta_0(V+\tau_0)=\mathrm{pr}_{[\lceil n/2\rceil,\lfloor n/2 \rfloor+2]}(d_{H_0}(V)+2)-(d_{H_0}(V)+2).
\]
Since $X(V)\le n$ by Lemma \ref{lnonnegdef}, $\Delta_{0}(V)=\max(0,\lceil n/2\rceil-d_{H_0}(V))$.

Therefore, we have
\[
\Delta_0(V+\tau_0)=\begin{cases}
 \Delta_0(V)-2 &   d_{H_0}(V) \le \lceil n/2\rceil-2,
 \\ \Delta_0(V)=0 &  \lceil n/2\rceil \le d_{H_0}(V) \le \lfloor n/2\rfloor, 
 \\ -1 &  d_{H_0}(V) = \lfloor n/2 \rfloor+1,
 \\ -2 &  d_{H_0}(V) = \lfloor n/2 \rfloor+2,
 \\ 0 &  d_{H_0}(V)=\lceil n/2\rceil -1, \mbox{ then } \Delta_{0}(V)=1.
 \end{cases}
\]

Thus, except when \(d_{H_0}(V)=\lfloor n/2\rfloor+1\) or \(d_{H_0}(V)=\lceil n/2\rceil-1\), the defect vector preserves parity. However, since \(\Delta_0(V)\) is a positive odd integer, only the cases \(d_{H_0}(V)\leq \lceil n/2\rceil-2\) and \(d_{H_0}(V)=\lceil n/2\rceil-1\) can occur. The latter case occurs precisely when \(\Delta_0(V)=1\).

By $\frac{\Delta_{0}(V)+1}{2}$ many times repeated application of $\tau_0$, we derive that
\[
\Delta(V+\frac{\Delta_0(V)+1}{2}\tau_0)= (0, \Delta_1(V), \Delta_2(V)).
\]
Consequently, we compute $x_0(V+\frac{\Delta_0(V)+1}{2}\tau_0)=\frac{\Delta_1(V)+\Delta_2(V)}{2}, x_2(V+\frac{\Delta_0+1}{2}\tau_0)=\frac{\Delta_2(V)-\Delta_1(V)}{2}$ and $x_2(V+\frac{\Delta_0(V)+1}{2}\tau_0)=\frac{\Delta_1(V)-\Delta_2(V)}{2}$. Since we want all the signs to be preserved, this only happens if and only if  $\Delta_1(V)=\Delta_2(V).$ If $\Delta_1(V)\neq\Delta_2(V)$, we add multiples of $\tau_1$ and $\tau_2$ so that the new values of $x_1$ and $x_2$ are zero. Therefore, we choose:
\begin{myeq}\label{vprime}
    V'= V+ \frac{\Delta_0(V)+1}{2}\tau_0 +\frac{\Delta_1(V)-\min(\Delta_1(V), \Delta_2(V))}{2}\tau_1+\frac{\Delta_2(V)-\min(\Delta_1(V), \Delta_2(V))}{2}\tau_2.
\end{myeq}

Hence, we derive that
\[
\Delta(V')=(0, \min(\Delta_1(V), \Delta_2(V)), \min(\Delta_1(V), \Delta_2(V)))\] 
and
\[(x_0(V'), x_1(V'), x_2(V'))= (\min(\Delta_1(V), \Delta_2(V)), 0, 0).
\]
 Thus $\sum_{\Delta(V')}$ is even and $x_i(V')\ge 0.$ 

Given $X(V) \le n$, we may observe that $X(V')\le n$ as well. To verify this, let us consider the fixed point dimensions $d_{H_i}(V')$ for $i \in [2]$. For simplicity, let us assume that $\Delta_2(V) \le \Delta_1(V)$. From \eqref{vprime}, we have:
\[
d_{H_0}(V')=d_{H_0}(V)+ \Delta_0(V)+1= \delta_0(V) +1,\quad d_{H_1}(V')=d_{H_1}(V)+\Delta_1(V)- \Delta_2(V), \quad d_{H_2}(V')=d_{H_2}(V).
\]

Since $X(V)\le n$, by Lemma \ref{lnonnegdef}
\begin{myeq}
\Delta_i(V)=\max(0,\,\lceil n/2\rceil-d_{H_i}(V))\ge 0 \mbox{ for } i\in [2].
\end{myeq}
 Since $\Delta_0(V)$ is odd, therefore $\Delta_0(V)\ge1$, thus $\Delta_0(V)=\lceil n/2\rceil -d_{H_0}(V).$ Hence $d_{H_0}(V')= \lceil n/2\rceil +1 \le \lfloor n/2\rfloor+2.$

 For $\Delta_1(V),$ if $\Delta_1(V) \le 1$, then $\Delta_1(V)=\lceil n/2\rceil -d_{H_1}(V)$, thus $d_{H_1}(V')= \lceil n/2\rceil -\Delta_{2}(V) \le \lfloor n/2\rfloor+2.$ For $\Delta_{1} (V)=0,$ then $0\le\Delta_{2}(V)\le\Delta_{1}(V)$ forces $\Delta_2(V)=0$. Hence, $d_{H_1}(V')=d_{H_1}(V).$

 Thus by Theorem \ref{main}, consider
 \[
 \widehat{V'}= V'+\min(\Delta_1(V), \Delta_2(V))(1-\alpha_0).
 \]
 Therefore, we have $P^n_n(X(V'))\simeq X(\widehat{V'}).$ Hence we have the cofibre sequence
\[
X(\widehat{V'})\xrightarrow{\ v_0\ }P^n_nX(V)\to P_n(\Sigma^{\widehat{V'}-\tau_0}C(v_0)),
\]
in which $X(\widehat{V'})$ is an $n$-slice. In particular, $P^n_nX(V)$ is not spherical.

(b) In this case, \(\Delta_0(V)\) is odd. Since \(X(V)\geq n\), Lemma~\ref{lnonnegdef} implies that
\[
\Delta_0(V)=\min\left(0,\left\lfloor n/2\right\rfloor+2-d_{H_0}(V)\right)\leq -1.
\]

Define 
\[
V'= V - \frac{|\Delta_0(V)|-1}{2}\tau_0- \frac{|\Delta_1(V)|-\min(|\Delta_1(V)|, |\Delta_2(V)|)}{2}\tau_1-\frac{|\Delta_2|-\min(|\Delta_1(V)|, |\Delta_2(V)|)}{2}\tau_2
\]
and
\[
\widehat{V'}= V'-\min(|\Delta_1(V)|, |\Delta_2(V)|)(1-\alpha_0). 
\]

Hence the result follows. \qedhere
\end{proof}

Using the above results, we completely recover \cite[Proposition~5.5]{GY20} as follows.

\begin{thm}{\label{GYThm}}
   (a) Let $n\ge 5$. Then $X(n) \ge n$, $\Delta(X(n))= (-k, -k, -k)$ and $x_i=-\frac{k}{2}$ where $i=0,1,2$ and $k = \lceil n/2\rceil-2.$ 
    \begin{enumerate}
        \item If $n=0,3 \pmod 4$,  then $P^n_nX(n)$ is spherical and $P^n_n X(n)\simeq X(n-\frac{k}{2}(3-\bar{\rho})).$
        \item If $n=1,2 \pmod 4,$ then $P^n_n X(n)$ is not spherical and there is cofibre sequence
        \[
         P^n(\Sigma^{\widehat{V'}-1}C(v_0)) \to P^n_n X(n) \xrightarrow{\ v_0\ }X\bigl(\widehat{V'}\bigr) 
        \]
        where 
        \begin{myeq}{\label{cons.V}}
        \widehat{V'}=(n-\frac{\lceil n/2\rceil-1}{2}-k)+(k-\frac{\lceil n/2\rceil-1}{2})\,\alpha_0+\frac{\lceil n/2\rceil-1}{2}\,\alpha_1+\frac{\lceil n/2\rceil-1}{2}\,\alpha_2.
        \end{myeq}

    Moreover, the $n$-slice is
    \[
    P^n_nX(n)\simeq S^{y\rho+(n-4y)} \wedge H\uf, \text{ where }y=\frac{\lceil n/2\rceil-1}{2}.
    \]
    \end{enumerate}

    (b) For $0\le n\le 4$, $X(n)$ is an $n$-slice.
\end{thm}
\begin{remark}
Although the construction of \(\widehat V'\) in~\eqref{cons.V} depends on the choice of \(\alpha_0\), replacing \(\alpha_0\) with any other \(\alpha_i\) does not change the resulting \(n\)-slice, by Lemma~\ref{kerequiv}.
\end{remark}
\begin{proof}[Proof of Proposition~\ref{GYThm}]
   (a) (1) It is straightforward to verify that \(\Delta(X(n))\) is even if and only if \(n\equiv 0,3\pmod 4\). Hence, by Theorem~\ref{main}, \(X(n)\) admits a spherical slice.

   (2) We first show that $P^n_nX(n) \simeq S^{\widehat{V'}}\wedge H\ker(\uF \twoheadrightarrow \phi^*_{H_0}\uF).$ 
    
    First consider a pair $(a,b)$ given by
    \[
(a,b)=
\begin{cases}\bigl(\frac{n+7}{4},\ \frac{n-5}{4}\bigr) & n\equiv1\pmod4,\\[4pt]
\bigl(\frac{n+10}{4},\ \frac{n-6}{4}\bigr) & n\equiv2\pmod4.
\end{cases}
\]
 Then using Corollary \ref{csphcvi}, we get a cofibre sequence   
    \[
\phi^*_{H_0}\bigl(S^{a+b\sigma}\wedge H_{C_2}\uF\bigr) \longrightarrow\Sigma^{\widehat{V'}-1}C(v_0)\longrightarrow \phi^*_{H_0}\bigl(S^{(a-1)+b\sigma}\wedge H_{C_2}\uF\bigr).
\]
 It follows from Proposition~\ref{C_2slices} and~\eqref{p:slices} that \(\phi^*_{H_0}\bigl(S^{(a-1)+b\sigma}\wedge H_{C_2}\uF\bigr)\) is a \(2\lfloor n/2\rfloor\)-slice, which is an \((n-1)\)-slice when \(n\) is odd and an \(n\)-slice when \(n\) is even. Similarly, \(\phi^*_{H_0}\bigl(S^{a+b\sigma}\wedge H_{C_2}\uF\bigr)\) is a \((2\lfloor n/2\rfloor+2)\)-slice. So $\phi^*_{H_0}\bigl(S^{a+b\sigma}\wedge H_{C_2}\uF\bigr) \ge n+1$.
Therefore
\[
P^n\bigl(\Sigma^{\widehat{V'}-1}C(v_0)\bigr)\simeq \phi^*_{H_0}\bigl(S^{(a-1)+b\sigma}\wedge H_{C_2}\uF\bigr) .
\]
Hence we obtained the following cofibre sequence
\[
         \phi^*_{H_0}\bigl(S^{(a-1)+b\sigma}\wedge H_{C_2}\uF\bigr) \to P^n_n X(n) \xrightarrow{\ v_0\ }X\bigl(\widehat{V'}\bigr)  \stackrel{\partial}{\to} \Sigma \phi^*_{H_0}\bigl(S^{(a-1)+b\sigma}\wedge H_{C_2}\uF\bigr).
        \]
The connecting map $\partial$ is classified by the homotopy group 
\[
\pi^{K_4}_{\widehat{V'}}(\Sigma \phi^*_{H_0}\bigl(S^{(a-1)+b\sigma}\wedge H_{C_2}\uF\bigr))\cong\pi^{C_2}_{(a+b\sigma)-1-(a-1)-b\sigma}H_{C_2}\uF=\pi^{C_2}_0H_{C_2}\uF\cong\F .
\]

Note that $\partial$ is not a null-homotopy.  On the contrary, assume that $\partial\sim \ast$. Then $\phi^*_{H_0}\bigl(S^{(a-1)+b\sigma}\wedge H_{C_2}\uF\bigr)$ is a retract of $P^n_nX(n).$ This is a contradiction as $P^n_n X(n) \ge n.$ Hence we compute that 
\[
P^n_nX(n) \simeq S^{\widehat{V'}}\wedge H\ker(\uF \twoheadrightarrow \phi^*_{H_0}\uF).
\]

Next using Lemma \ref{kerequiv}, we have the equivalence $H\ker(\uF \twoheadrightarrow \phi^*_{H_0}\uF) \simeq \Sigma^{\alpha_0-1}\wedge H\uf$. This yields
    \[
    P^n_nX(n) \simeq S^{\widehat{V'}+\alpha_0-1}\wedge H\uf=S^{y\rho+(n-4y)} \wedge H\uf,
    \]
    where $y=\frac{\lceil n/2\rceil-1}{2}.$ Hence the result follows.

    (b) The result follows directly from Corollary~\ref{pk4range}.
\end{proof}

The following is an example of a spherical slice.
\begin{example}
    Consider the following spectrum
    \[
    X(2\aO)=S^{2\aO}\wedge H\uF.
    \]
    Here $|V|=2$ and
$(d_{H_0}(2\alpha_0),d_{H_1}(2\alpha_0),d_{H_2}(2\alpha_0),d_{K_4}(2\alpha_0))=(2,0,0,0)$, so $X(2\aO)\in[n(V),k(V)]=[0,2]$. The target dimensions are
$(\delta_{H_0}(2\alpha_0),\delta_{H_1}(2\alpha_0),\delta_{H_2}(2\alpha_0))=(2,1,1)$ and the defect vector is $\Delta(2\alpha_0)=(0,1,1)$. Its sum is even,
and $x_0(2\alpha_0)=\tfrac12(\Delta_1(2\alpha_0)+\Delta_2(2\alpha_0)-\Delta_0(2\alpha_0))=1$, $x_1(2\alpha_0)=x_2(2\alpha_0)=0$. So the $2$-slice should be $X(U)$
with $U=2\aO+(1-\aO)=1+\aO$, reached by one raising move along $u_{\aO}$. Thus, this gives the slice tower
\[
\xymatrix{
X(1+\aO)\simeq P^2_2X(2\aO)\ar[r]^-{u_{\aO}} & X(2\aO)\simeq P^2X(2\aO)\ar[d]\\
 & P^1X(2\aO)\simeq P^0_0X(2\aO)\simeq H\ug_{(0)}.
}
\]
\end{example}

 Next, we have an example that is not a spherical slice.

\begin{example}\label{botx5}
    Consider the $K_4$-spectrum $X(5)=\Sigma^5H\uF$ which lies in $[5,8]$ by Corollary~\ref{bound0}. Note that
    $$\Delta(5)=(\Delta_{0}(5),\Delta_1(5),\Delta_2(5))=(-1,-1,-1),\quad \sum\nolimits_{\Delta(5)}=-3.$$
Hence we have $$x_0(5)= x_1(5)=x_2(5)=-\tfrac12 .$$
By Theorem \ref{mainodd}, take $$V'=V-\tau_0=4-\alpha_0+\alpha_1+\alpha_2 .$$ 

 Hence $\widehat{V'}= V'-(1-\alpha_0)=3+\alpha_1+\alpha_2.$ Thus $P^5_5(X(V'))= S^{3+\alpha_1+\alpha_2}\wedge H\uF.$ Therefore we have a cofibre sequence from \eqref{odd2} 
\[
P^5\bigl(\Sigma^{2+\aI+\aII}C(v_0)\bigr)\longrightarrow P^5_5X(5)\longrightarrow S^{3+\aI+\aII}\wedge H\uF.
\]

For $\Sigma^{2+\aI+\aII}C(v_0)$, using Corollary \ref{csphcvi}, we get a cofibre sequence
\[
\Sigma^3 H \phi^*_{H_0}\uF \to \Sigma^{2+\aI+\aII}C(v_0) \to  \Sigma^2 H \phi^*_{H_0}\uF.
\]
 Here $\Sigma^2 H \phi^*_{H_0}\uF$ is a $4$-slice and $\Sigma^3 H \phi^*_{H_0}\uF$ is a $6$-slice. Hence $P^5\bigl(\Sigma^{2+\aI+\aII}C(v_0)\bigr)\simeq \Sigma^2 H \phi^*_{H_0}\uF$.
 Hence the $5$-slice of $X(5)$ sits inside the cofibre sequence:
 $$\Sigma^2H\phi^*_{H_0}\uF\longrightarrow P^5_5X(5)\longrightarrow S^{3+\alpha_1+\alpha_2}\wedge H\uF .$$
 Now we compute the connecting map, which is classified by the homotopy class $[S^{3+\alpha_1+\alpha_2}\wedge H\uF, \Sigma^3 H\phi^*_{H_0}\uF]\cong \pi_{\alpha_1+\alpha_2}^{K_4}(H\phi^*_{H_0}\uF)\cong \F.$ If this cofibre sequence splits, then $\Sigma^2H\phi^*_{H_0}\uF$ is a retract of $P^5_5X(5)$, which gives $\Sigma^2H\phi^*_{H_0}\uF \ge 5$. This is a contradiction. Hence the $5$-slice is
 \[
 P^5_5 X(5)\simeq S^{3+\alpha_1+\alpha_2}\wedge H\ker(\uF \twoheadrightarrow \phi^*_{H_0}\uF )
 \]
 where the Lewis diagram of the Mackey functor $\ker(\uF \twoheadrightarrow \phi^*_{H_0}\uF )$  is given as $\underline{K}_{H_0}$ in~\eqref{makK}.
\end{example}

\section{Appendix. Mackey functors}\label{sapp}In this section, we explicitly describe all the Mackey functors used throughout this article. We represent these Mackey functors using Lewis diagrams~\cite{L88}. For convenience, all maps omitted from the below diagrams are understood to be zero maps.

\begin{mysubsection}{$C_2$-Mackey functors}
~Let $\sigma$ denote the sign representation of $C_2$.

\[
\begin{array}{|c@{\hspace{0.8cm}}|c@{\hspace{0.8cm}}|c@{\hspace{0.8cm}}|c|}
\hline
\xymatrix@R=1.1cm{
\F\ar@/_.5pc/[d]^{1}\\
\F
}
&
\xymatrix@R=1.1cm{
\F\\
\F\ar@/_.5pc/[u]^{1}
}
&
\xymatrix@R=1.1cm{
\F\\
0
}
&
\xymatrix@R=1.1cm{
0\\
\F
}
\\[0.4cm]
\hline
\uF
&
\uF^{\ast}
&
\underline g_{C_2}
&
\underline f_{C_2}
\\
\hline
\end{array}
\]
These are related by
\[
\Sigma^{2-2\sigma}H_{C_2}\uF\simeq H_{C_2}\uF^{\ast},\qquad
\Sigma^{1-\sigma}H_{C_2}\uF\simeq H_{C_2}\underline f_{C_2},\qquad
\Sigma^{\sigma}H_{C_2}\underline g_{C_2}\simeq H_{C_2}\underline g_{C_2}.
\]
In particular $\Sigma^{\rho_{C_2}}H_{C_2}\underline g_{C_2}\simeq\Sigma H_{C_2}\underline g_{C_2}$.
\end{mysubsection}

\begin{mysubsection}{$K_4$-Mackey functors}
~Write $H_i=\ker\alpha_i$ for $i=0,1,2$. In each diagram, the middle row lists the values at $K_4/H_0$, $K_4/H_1$ and $K_4/H_2$,
from left to right, and the bottom entry is the value at $K_4/e$.

\begin{example}
    \[
\begin{array}{|c|c|}
\hline
\xymatrix@C=3em@R=3em{
 & \F
   \ar@/_0.5pc/[ld]_{1}
   \ar@/^0.5pc/[d]^{1}
   \ar@/^0.5pc/[rd]^{1}
 & \\
\F
   \ar@/_0.5pc/[rd]_{1}
 & \F
   \ar@/^0.5pc/[d]^{1}
 & \F
   \ar@/^0.5pc/[ld]^{1}
 \\
 & \F
 &
}
&
\xymatrix@C=3em@R=3em{
 & \F
 & \\
\F
   \ar@/^0.5pc/[ru]_{1}
 & \F
   \ar@/^0.5pc/[u]_{1}
 & \F
   \ar@/_0.5pc/[lu]^{1}
 \\
 & \F
   \ar@/^0.5pc/[lu]_{1}
   \ar@/^0.5pc/[u]_{1}
   \ar@/_0.5pc/[ru]^{1}
 &
}
\\[0.5cm]
\hline
\uF & \uF^\ast
\\
\hline
\end{array}
\]
They satisfy $\Sigma^{3-\bar{\rho}}H\uF\simeq H\uF^*.$
\end{example}

\begin{example} Consider the following $K_4$-Mackey functors

\[
\begin{array}{|c|c|}
\hline
\xymatrix@C=2.5em@R=2.5em{
 & \F & \\
0 & 0 & 0 \\
 & 0 &
}
&
\xymatrix@C=3em@R=3em{
 & 0 & \\
0 & 0 & 0 \\
 & \F &
}
\\[0.5cm]
\hline
\ug & \uf
\\
\hline
\end{array}
\]
For every $V\in RO(K_4)$ we have $\Sigma^VH\ug\simeq\Sigma^{\dim V^{K_4}}H\ug$.
\end{example}

\begin{example}These Mackey functor obtained as the pullback Mackey functor of $C_2$-Mackey functor $\uF$.
\[
\begin{array}{|c|c|c|}
\hline
\xymatrix@C=0.7cm@R=0.8cm{
  & \F \ar@/_.5pc/[dl]^{1}\\
  \F & 0 & 0 \\
  & 0
}
&
\xymatrix@C=0.7cm@R=1cm{
  & \F \ar@/_.5pc/[d]^{1}\\
  0 & \F & 0 \\
  & 0
}
&
\xymatrix@C=0.7cm@R=1cm{
  & \F\ar@/^.5pc/[dr]^{1} \\
  0 & 0 & \F \\
  & 0
}
\\[0.6cm]
\hline
\phi^{\ast}_{H_0}(\uF)
&
\phi^{\ast}_{H_1}(\uF)
&
\phi^{\ast}_{H_2}(\uF)
\\
\hline
\end{array}
\]
\end{example}

\begin{example}These Mackey functor obtained as the pullback Mackey functor of $C_2$-Mackey functor $\uF^{\ast}$.
\[
\begin{array}{|c|c|c|}
\hline
\xymatrix@C=0.7cm@R=1cm{
  & \F \\
  \F \ar@/^.5pc/[ur]^{1} & 0 & 0 \\
  & 0
}
&
\xymatrix@C=0.7cm@R=1cm{
  & \F \\
  0 & \F \ar@/^.5pc/[u]^{1} & 0 \\
  & 0
}
&
\xymatrix@C=0.7cm@R=1cm{
  & \F \\
  0 & 0 & \F\ar@/_.5pc/[ul]^{1} \\
  & 0
}
\\[0.6cm]
\hline
\phi^{\ast}_{H_0}(\uF^{\ast})
&
\phi^{\ast}_{H_1}(\uF^{\ast})
&
\phi^{\ast}_{H_2}(\uF^{\ast})
\\
\hline
\end{array}
\]
\end{example}

\begin{example} The following Mackey functors are obtained as the pullback Mackey of $\underline{f}_{C_2}$ along the order $2$ subgroups. 
\[
\begin{array}{|c|c|c|}
\hline
\xymatrix@C=0.7cm@R=1cm{
  & 0 \\
  \F & 0 & 0 \\
  & 0
}
&
\xymatrix@C=0.7cm@R=1cm{
  & 0 \\
  0 & \F & 0 \\
  & 0
}
&
\xymatrix@C=0.7cm@R=1cm{
  & 0 \\
  0 & 0 & \F \\
  & 0
}
\\[0.6cm]
\hline
\ug_{H_0}:=\phi^{\ast}_{H_0}(\underline{f}_{C_2})
&
\ug_{H_1}:=\phi^{\ast}_{H_1}(\underline{f}_{C_2})
&
\ug_{H_2}:=\phi^{\ast}_{H_2}(\underline{f}_{C_2})
\\
\hline
\end{array}
\]
\end{example}

\begin{example}{\label{MacN}} The following Mackey functors appear  in~\eqref{dcvi}.
\[
\begin{array}{|c|c|c|}
\hline
\xymatrix@C=0.9cm@R=1.1cm{
  & \F \ar@/_.5pc/[dl]^{1} \\
  \F
  & \F\ar@/_.5pc/[d]^{1} \ar@/_.5pc/[u]^{1}
  & \F\ar@/^.5pc/[dl]^{1}
  \ar@/_.5pc/[ul]^{1}\\
  & \F\ar@/^.5pc/[ul]^{1}
}
&
\xymatrix@C=0.9cm@R=1.1cm{
  & \F \ar@/_.5pc/[d]^{1} \\
  \F \ar@/_.5pc/[dr]^{1}
  \ar@/^.5pc/[ur]^{1}
  & \F
  & \F\ar@/^.5pc/[dl]^{1}
  \ar@/_.5pc/[ul]^{1}\\
  & \F\ar@/_.5pc/[u]^{1}
}
&
\xymatrix@C=0.9cm@R=1.1cm{
  & \F \ar@/^.5pc/[dr]^{1} \\
  \F\ar@/_.5pc/[dr]^{1}
  \ar@/^.5pc/[ur]^{1}
  & \F\ar@/_.5pc/[d]^{1}
  \ar@/_.5pc/[u]^{1}
  & \F \\
  & \F\ar@/_.5pc/[ur]^{1}
}
\\[0.5cm]
\hline
\underline{N}_{H_0}
&
\underline{N}_{H_1}
&
\underline{N}_{H_2}
\\
\hline
\end{array}
\]
\medskip
They satisfy the equivalence $H\uN_{H_i}\simeq \Sigma^{\tau_i}H\uF.$
\end{example}

\begin{example}{\label{Phimac}}These Mackey functors are obtained by inflating the \(C_2\)-Mackey functor \(\ug_{C_2}\) from Corollary~\ref{tformula}. They also arise in Proposition~\ref{pCform} as the cofibre of the map induced by multiplication by the Thom class. 
\[
\begin{array}{|c|c|c|}
\hline
\xymatrix@C=0.7cm@R=1cm{
  & \F\ar@/_.5pc/[dr]^{1}\ar@/_.5pc/[d]^{1} \\
  0 & \F & \F \\
  & 0
}
&
\xymatrix@C=0.7cm@R=1cm{
  & \F \ar@/_.5pc/[dl]^{1}\ar@/_.5pc/[dr]^{1}\\
  \F & 0 & \F \\
  & 0
}
&
\xymatrix@C=0.7cm@R=1cm{
  & \F\ar@/_.5pc/[dl]^{1}\ar@/_.5pc/[d]^{1} \\
  \F & \F & 0 \\
  & 0
}
\\[0.5cm]
\hline
\ug_{(0)}:=\Phi^*_{H_0}\ug_{C_2}
&
\ug_{(1)}:=\Phi^*_{H_1}\ug_{C_2}
&
\ug_{(2)}:=\Phi^*_{H_2}\ug_{C_2}
\\
\hline
\end{array}
\]
\end{example}
\end{mysubsection}


\begin{thebibliography}{99}

\bibitem{BG21}
{\scshape S.~Basu, S.~Ghosh,}
\emph{Bredon cohomology of finite dimensional $C_p$-spaces,}
Homology Homotopy Appl.\ {\bf 23} (2021), no.~2, 33--57.

\bibitem{BHLSZ22}
{\scshape A.~Beaudry, M.~A.~Hill, T.~Lawson, X.~D.~Shi, M.~Zeng,}
\emph{On the slice spectral sequence for quotients of norms of Real bordism,}
J.\ Topol.\ {\bf 18} (2025), no.~1, Paper No.~e70015.

\bibitem{BHSZ21}
{\scshape A.~Beaudry, M.~A.~Hill, X.~D.~Shi, M.~Zeng,}
\emph{Models of Lubin--Tate spectra via Real bordism theory,}
Adv.\ Math.\ {\bf 392} (2021), Paper No.~108020.

\bibitem{Gho22}
{\scshape S.~Ghosh}, \emph{Slice Filtration of certain $C_{pq}$ Spectra,} New York J.\ Math.\ {\bf 28} (2022), 1399–1418.

\bibitem{GK26}
{\scshape S.~Ghosh, A.~Kumar,}
\emph{$RO(C_p\times C_p)$-graded cohomology of universal spaces and the coefficient ring,}
preprint, arXiv:2603.10626 (2026).

\bibitem{GM95}
{\scshape J.~P.~C.~Greenlees, J.~P.~May,}
\emph{Equivariant stable homotopy theory,}
in Handbook of Algebraic Topology, North-Holland, Amsterdam, 1995, 277--323.

\bibitem{GKM25}
{\scshape B.~J.~Guillou, J.~Keyes, D.~Mehrle,}
\emph{The $RO(K)$-graded homotopy of Klein-four normed Mackey functors,}
Glasg.\ Math.\ J.\ (2026), doi:10.1017/S0017089526101001.

\bibitem{GKM26}
{\scshape B.~J.~Guillou, J.~Keyes, D.~Mehrle,}
\emph{The zero slice of quaternionic Real bordism,}
preprint, arXiv:2605.01174 (2026).

\bibitem{GS25}
{\scshape B.~J.~Guillou, C.~Slone,}
\emph{The slices of quaternionic Eilenberg--Mac Lane spectra,}
Algebr.\ Geom.\ Topol.\ {\bf 25} (2025), no.~6, 3341--3384.

\bibitem{GY20}
{\scshape B.~J.~Guillou, C.~Yarnall,}
\emph{The Klein four slices of $\Sigma^nH\underline{\mathbb F}_2$,}
Math.\ Z.\ {\bf 295} (2020), no.~3--4, 1405--1441.

\bibitem{Hil12}
{\scshape M.~A.~Hill,}
\emph{The equivariant slice filtration: a primer,}
Homology Homotopy Appl.\ {\bf 14} (2012), no.~2, 143--166.

\bibitem{HHR16}
{\scshape M.~A.~Hill, M.~J.~Hopkins, D.~C.~Ravenel,}
\emph{On the nonexistence of elements of Kervaire invariant one,}
Ann.\ of Math.\ (2) {\bf 184} (2016), no.~1, 1--262.

\bibitem{HHR17b}
{\scshape M.~A.~Hill, M.~J.~Hopkins, D.~C.~Ravenel,}
\emph{The slice spectral sequence for certain $RO(C_{p^n})$-graded suspensions of $H\underline{\mathbb Z}$,}
Bol.\ Soc.\ Mat.\ Mex.\ (3) {\bf 23} (2017), no.~1, 289--317.

\bibitem{HSWX23}
{\scshape M.~A.~Hill, X.~D.~Shi, G.~Wang, Z.~Xu,}
\emph{The slice spectral sequence of a $C_4$-equivariant height-$4$ Lubin--Tate theory,}
Mem.\ Amer.\ Math.\ Soc.\ {\bf 288} (2023), no.~1429.

\bibitem{HY18}
{\scshape M.~A.~Hill, C.~Yarnall,}
\emph{A new formulation of the equivariant slice filtration with applications to $C_p$-slices,}
Proc.\ Amer.\ Math.\ Soc.\ {\bf 146} (2018), no.~8, 3605--3614.


\bibitem{HK17}
{\scshape J.~Holler, I.~Kriz.} 
\emph{``On $RO(G)$–graded equivariant ``ordinary” cohomology where $G$ is
a power of $\mathbb Z/2$."} \emph{Algebr. Geom. Topol.} {\bf 17.2} (2017): 741-763.

\bibitem{K21}
{\scshape S.~Kriz,}
\emph{Notes on equivariant homology with constant coefficients,}
Pacific J.\ Math.\ {\bf 309} (2021), no.~2, 381--399.

\bibitem{L88}
{\scshape L.~G.~Lewis,}
\emph{The $RO(G)$-graded equivariant ordinary cohomology of complex projective spaces with linear $\mathbb Z/p$-actions,}
in Algebraic Topology and Transformation Groups (G\"ottingen, 1987), Lecture Notes in Math.\ {\bf 1361}, Springer, Berlin,
1988, 53--122.

\bibitem{LSWX19}
{\scshape G.~Li, X.~D.~Shi, G.~Wang, Z.~Xu,}
\emph{Hurewicz images of real bordism theory and real Johnson--Wilson theories,}
Adv.\ Math.\ {\bf 342} (2019), 67--115.

\bibitem{May96}
{\scshape J.~P. May}, {\em Equivariant homotopy and cohomology theory}, vol.~91 of
  CBMS Regional Conference Series in Mathematics, Published for the Conference
  Board of the Mathematical Sciences, Washington, DC; by the American
  Mathematical Society, Providence, RI, 1996.
\newblock With contributions by M. Cole, G. Comeza\~na, S. Costenoble, A. D.
  Elmendorf, J. P. C. Greenlees, L. G. Lewis, Jr., R. J. Piacenza, G.
  Triantafillou, and S. Waner.

\bibitem{MM02}
{\scshape M.~A.~Mandell, J.~P.~May,}
\emph{Equivariant orthogonal spectra and $S$-modules,}
Mem.\ Amer.\ Math.\ Soc.\ {\bf 159} (2002), no.~755.

\bibitem{MSZ20}
{\scshape L.~Meier, X.~D.~Shi, M.~Zeng,}
\emph{The localized slice spectral sequence, norms of Real bordism, and the Segal conjecture,}
Adv.\ Math.\ {\bf 412} (2023), Paper No.~108804.

\bibitem{SY25}
{\scshape I.~Sikora, G.~Yan,}
\emph{On the structure of the $RO(G)$-graded homotopy of $H\underline M$ for cyclic $p$-groups,}
Algebr.\ Geom.\ Topol.\ {\bf 25} (2025), 2933--2980.

\bibitem{Slo22}
{\scshape C.~Slone,}
\emph{Klein four $2$-slices and the slices of $\Sigma^{\pm n}H\underline{\mathbb Z}$,}
Math.\ Z.\ {\bf 301} (2022), 3895--3938.

\bibitem{U}
{\scshape J.~R.~Ullman,}
\emph{On the regular slice spectral sequence,}
Ph.D.\ thesis, Massachusetts Institute of Technology, 2013.

\bibitem{Yar17}
{\scshape C.~Yarnall,}
\emph{The slices of $S^n\wedge H\underline{\mathbb Z}$ for cyclic $p$-groups,}
Homology Homotopy Appl.\ {\bf 19} (2017), no.~1, 1--22.

\bibitem{Zen17}
{\scshape M.~Zeng,}
\emph{Equivariant Eilenberg--Mac Lane spectra in cyclic $p$-groups,}
preprint, arXiv:1710.01769 (2017).

\end{thebibliography}
\end{document}